# A FRAMEWORK FOR COUPLED DEFORMATION-DIFFUSION ANALYSIS WITH APPLICATION TO DEGRADATION/HEALING

M. K. MUDUNURU AND K. B. NAKSHATRALA

ABSTRACT. This paper deals with the formulation and numerical implementation of a fully coupled continuum model for deformation-diffusion in linearized elastic solids. The mathematical model takes into account the effect of the deformation on the diffusion process, and the affect of the transport of an inert chemical species on the deformation of the solid. We then present a robust computational framework for solving the proposed mathematical model, which consists of coupled non-linear partial differential equations. It should be noted that many popular numerical formulations may produce *unphysical* negative values for the concentration, particularly, when the diffusion process is anisotropic. The violation of the non-negative constraint by these numerical formulations is not mere numerical noise. In the proposed computational framework we employ a novel numerical formulation that will ensure that the concentration of the diffusant be always non-negative, which is one of the main contributions of this paper. Representative numerical examples are presented to show the robustness, convergence, and performance of the proposed computational framework. Another contribution of this paper is to systematically study the affect of transport of the diffusant on the deformation of the solid and vice-versa, and their implication in modeling degradation/healing of materials. We show that the coupled response is both qualitatively and quantitatively different from the uncoupled response.

## 1. INTRODUCTION

In this paper we present a continuum mathematical model and a computational framework for degradation/healing in elastic solids due to the presence of a solute. We shall neglect chemical reactions as well as thermal effects. We shall also assume that the strains in the solid are small. The deformation is coupled with the diffusion process, and the diffusion process is in turn coupled with the deformation of the solid. The chosen problem belongs to a broader class of problems, namely, coupled deformation-diffusion problems.

Coupled deformation-diffusion problems arise in many civil engineering, material science, and polymer science applications. An important example is degradation/healing of materials and structures. Many man-made and natural materials degrade/heal due to environmental conditions, and structural components and superstructures are constantly exposed to adverse conditions. The fate of the transport of a diffusant will in turn depend on the deformation of the solid. Some other examples of coupled deformation-diffusion problems are moisture damage in cementitious materials [30, 31] and asphalt [11, 26], hydrogen embrittlement [51, 57], curing of ceramics [35], and swelling of polymers and composites [58, 55].

The governing equations of a deformation-diffusion problem are coupled non-linear partial differential equations. That is, the mechanical variables (strains and stresses) and diffusive variables (concentration and mass flux) are coupled through constitutive relations. It is (in general) not possible to obtain analytical solutions for these kinds of problems, and one has to resort to numerical techniques to solve practical problems. Despite the importance of coupled deformation-diffusion problems, there is no *robust* and *reliable* computational framework for solving such coupled problems. In particular, the existing numerical studies have one or more of the following limitations:

- considered academic and unrealistic problems like infinite slabs and infinite cylinders,
- did not consider anisotropic diffusion and/or assumed the medium to be homogeneous, or
- did not consider the fact that conventional numerical formulations and finite element packages produce unphysical negative values for the concentration in solving diffusion-type equations.

1.1. **Maximum principles and non-negative solutions for diffusion-type equations.** Predictive numerical simulations require accurate and reliable discretization methods. The resulting discrete systems must inherit or mimic fundamental properties of continuous systems. Maximum principles form an important set of properties for diffusion-type equations as these maximum principles have mathematical implications and physical consequences. In the study on partial differential equations, maximum principles are often used in existence theorems, and in obtaining point-wise estimates. For further details on (continuous) maximum principles refer to [40, 15, 18, 17, 47, 20].

A direct consequence of maximum principles for diffusion-type is the non-negativity of the solution (under appropriate conditions on the source and boundary conditions). Physical quantities like concentration of a diffusant should be non-negative by their nature and their approximations should also be non-negative as well. The question to ask is whether a chosen numerical formulation satisfies these maximum principles and meet the non-negative constraint. The discrete version of maximum principles is commonly referred as *discrete maximum principles*.

Many existing numerical formulations and packages do not satisfy the maximum principles. They may produce negative values for the primary variables in diffusion-type equations (that is, negative values for the concentration and temperature). *It should be emphasized that the violation is not mere numerical noise, that is, the violations will be much larger than the machine precision and cannot be neglected.*

For example, in Figure 1, we have shown that the contours of concentration obtained using Abaqus [1] for pure diffusion. The uncoupled problem is similar to one considered in Section

4. As one can see from the figure, significant part of the domain has negative solution. The minimum value of the concentration is approximately $-0.0832$, which is 4.16% off the range of possible values (which is between 0 and 2). In a subsequent section we show that the classical single-field Galerkin formulation produces negative values for the concentration for the same test problem (see Figure 15), and produces (qualitatively and quantitatively) wrong results for a coupled deformation-diffusion problem. Furthermore, Nakshatrala and Valocchi [45] have shown that the lowest order Raviart-Thomas [48] and variational multiscale [36, 44] mixed formulations violate discrete maximum principle and produce negative solutions for pure diffusion equation. Nagarajan and Nakshatrala [42] have shown that the single-field formulation violates discrete maximum principles for diffusion with decay.

1.2. **Our approach.** Herein, we will consider realistic problems, allow the medium to be inhomogeneous, consider anisotropic diffusion, and develop a robust computational framework that will always produce physically meaningful non-negative values for the concentration of the diffusant on general computational grids. The computational framework will consist of a non-negative formulation for diffusion equation, a single-field formulation for the deformation problem, and a staggered coupling algorithm.

We employ a staggered coupling technique (also known as partitioned solution approach) to couple individual analyses to obtain the coupled response. The non-negative formulation for diffusion is developed by extending the Galerkin formulation using convex programming, which produces (physical) non-negative solutions for the concentration even on general computational grids with low-order finite elements (linear three-node triangular, bilinear four-node quadrilateral, linear four-node tetrahedron, and tri-linear eight-node brick elements). The proposed non-negative formulation being applicable only to low-order finite elements is not a limitation as low-order finite elements remain quite popular in the solution of practical problems. (This is particularly true for large-scale simulations with complex geometries because of the inherent simplicity of low-order elements. Another reason is that adaptive mesh-generation techniques are simpler and tend to perform better with low-order finite elements.)

1.3. **Main contributions of this paper.** Some of the main contributions of this paper are as follows:

(1) We presented a mathematical model for diffusion of an inert chemical species in an elastic deformable solid that takes into account the effect of diffusant on the deformation. The model is *truly coupled* in the sense that the deformation of the solid will be affected by the diffusion process, and the diffusion process is in turn affected by the deformation of the solid. Such a coupled deformation-diffusion model is suitable to study degradation/healing in elastic solids. We restricted our model to steady-state.

(2) We presented a robust computational framework for performing deformation–diffusion analysis. The framework includes a solver for deformation, a non-negative solver for tensorial-diffusion, and a coupling algorithm to couple the individual deformation and diffusion analyses. The coupling algorithm is a staggered algorithm in which deformation and diffusion sub-problems are solved in an iterative fashion until convergence. *An important aspect of the proposed framework is that it employs a novel numerical formulation that ensures non-negative solutions for the concentration of the diffusant.*

(3) Using the proposed computational framework, we solved some realistic finite domain problems on general computational meshes. Also, we systematically studied the effect of the concentration of the diffusant on the deformation of the solid and vice-versa, and their implications on degradation/healing of materials and structures.

1.4. **An outline of this paper and symbolic notation.** The remainder of this paper is organized as follows. Section 2 presents a mathematical model for degradation/healing of a deformable elastic solid. The mathematical model will require performing coupled deformation-diffusion analysis. In Section 3 we present a fully coupled computational framework for deformation-diffusion analysis. The proposed computational framework will contain a non-negative formulation for tensorial-diffusion equation on general computational grids (which will always produce physically meaningful non–negative values for the concentration), a numerical solver for deformation, and a staggered coupling algorithm for coupling individual diffusion and deformation numerical solvers. In Section 4, representative numerical examples will be presented to illustrate the performance of the proposed coupled deformation-diffusion computational framework. Conclusions are drawn in Section 5.

The symbolic notation adopted in this paper is as follows. We shall make a distinction between vectors in the continuum and finite element settings. Similarly, we make a distinction between second-order tensors in the continuum setting versus matrices in the context of the finite element method. The continuum vectors are denoted by lower case boldface normal letters, and the second-order tensors will be denoted using upper case boldface normal letters (for example, vector **u** and second-order tensor **T**). In the finite element context, we shall denote the vectors using lower case boldface italic letters, and the matrices are denoted using upper case boldface italic letters (for example, vector $\boldsymbol{v}$ and matrix $\boldsymbol{K}$). Other notational conventions adopted in this paper are introduced as needed.

## 2. A MATHEMATICAL MODEL FOR COUPLED DEFORMATION-DIFFUSION

Consider an *inert* (chemical or biological) species being diffused through a deformable elastic solid. We now present a simple mathematical model for such a process. Let $\Omega \subset \mathbb{R}^{nd}$ be a bounded open domain, where "$nd$" denotes the number of spatial dimensions. The boundary $\partial\Omega$ is assumed to be piecewise smooth. Mathematically, $\partial\Omega = \bar{\Omega} - \Omega$, where $\bar{\Omega}$ is the set closure of $\Omega$. A spatial point in $\bar{\Omega}$ is denoted by $\mathbf{x}$. The gradient and divergence operators with respect to $\mathbf{x}$ are denoted by $\mathrm{grad}[\cdot]$ and $\mathrm{div}[\cdot]$, respectively. The unit outward normal to the boundary is denoted by $\mathbf{n}(\mathbf{x})$. We shall denote the displacement of the solid by $\mathbf{u}(\mathbf{x})$, and the concentration of the diffusant by $c(\mathbf{x})$. It is important to note that the concentration is a non-negative quantity, and a robust numerical solver should meet the non-negative constraint (which is not the case with many popular numerical schemes).

For the deformation problem, the boundary is divided into two complementary parts: $\Gamma_u^{\mathrm{D}}$ on which the displacement vector is prescribed, and $\Gamma_u^{\mathrm{N}}$ on which the traction vector is prescribed. For the diffusion problem, the boundary is divided into $\Gamma_c^{\mathrm{D}}$ on which concentration is prescribed, and $\Gamma_c^{\mathrm{N}}$ on which the flux is prescribed. For well-posedness, we require that $\Gamma_u^{\mathrm{D}} \cap \Gamma_u^{\mathrm{N}} = \emptyset$, $\Gamma_u^{\mathrm{D}} \cup \Gamma_u^{\mathrm{N}} = \partial\Omega$, $\Gamma_c^{\mathrm{D}} \cap \Gamma_c^{\mathrm{N}} = \emptyset$, and $\Gamma_c^{\mathrm{D}} \cup \Gamma_c^{\mathrm{N}} = \partial\Omega$. In addition, we assume that $\mathrm{meas}(\Gamma_u^{\mathrm{D}}) > 0$ and $\mathrm{meas}(\Gamma_c^{\mathrm{D}}) > 0$ for uniqueness of the solution.

### 2.1. A model for diffusion-dependent deformation.

The solid is modeled using linearized elasticity but the material parameters are allowed to depend on the concentration. The linearized strain is defined through

$$\mathbf{E}_l := \frac{1}{2}\left(\mathrm{grad}[\mathbf{u}] + \mathrm{grad}[\mathbf{u}]^T\right) \tag{1}$$

For a given concentration $c(\mathbf{x})$, the stress-strain relationship will be modeled as follows:

$$\mathbf{T}_c(\mathbf{u}, \mathbf{x}) = \lambda(\mathbf{x}, c)\mathrm{tr}[\mathbf{E}_l]\mathbf{I} + 2\mu(\mathbf{x}, c)\mathbf{E}_l \tag{2}$$

where $\mathbf{T}_c$ is the Cauchy stress, and $\lambda$ and $\mu$ are the Lamé parameters but now can depend both on the concentration and position. A simple model for the Lamé parameters to account for degradation/healing of the material due to the presence of a diffusant can be taken as follows:

$$\lambda(\mathbf{x}, c) = \lambda_0(\mathbf{x}) + \lambda_1(\mathbf{x})\frac{c(\mathbf{x})}{c_{\mathrm{ref}}} \tag{3a}$$

$$\mu(\mathbf{x}, c) = \mu_0(\mathbf{x}) + \mu_1(\mathbf{x})\frac{c(\mathbf{x})}{c_{\mathrm{ref}}} \tag{3b}$$

where $c_{\mathrm{ref}}$ is the reference concentration (which depends on the problem); $\lambda_0$ and $\mu_0$ are the Lamé parameters for the virgin material (that is, in the absence of the diffusant); and $\lambda_1(\mathbf{x})$ and $\mu_1(\mathbf{x})$ are the weights to account for the effect of concentration on the Lamé parameters.

The material parameters $\lambda_1$ and $\mu_1$ can be individually positive (which means that the material is healing), negative (which means that the material is degrading), or zero (which means that the material is unaffected by the presence of the diffusant). However, it is assumed that the parameters ($\lambda_0(\mathbf{x})$, $\lambda_1(\mathbf{x})$, $\mu_0(\mathbf{x})$ and $\mu_1(\mathbf{x})$) and concentration are such that we have bulk modulus $\lambda(\mathbf{x}, c) + \frac{2}{3}\mu(\mathbf{x}, c) > 0$ and shear modulus $\mu(\mathbf{x}, c) > 0$.

The above model given by equation (3) is similar to the concept of macroscopic damage variables, which has been introduced by Kachanov [24] to model material damage. The basic idea behind the concept of macroscopic damage variables is to quantify damage using internal variable(s). This concept has now been widely employed in numerous other works on damage (for example, see the review articles by Chaboche [8, 9]). In the above model, the concentration of the diffusant can be thought as a macroscopic damage/healing variable. The most common criticism about using internal variables is that these variables (in many cases) cannot be measured using physical experiments. However, in the above model, the dependence of the Lamé parameters on the concentration can be measured indirectly by non-destructive testing methods.

It should be noted that Weitsman [56], Kringos *et al.* [29, 28], and Muliana *et al.* [41] have considered material degradation as a function of concentration of the diffusant. However, the models developed in these works are not fully coupled. That is, the material properties of the solid are dependent on the concentration but the diffusivity is unaffected by the deformation of the solid. In Reference [25], a fully coupled model is proposed but a specific boundary value problem (torsion of a cylindrical annulus undergoing degradation) is solved using the semi-inverse method. They did not present a computational framework, and also their mathematical model is different from the one present in this paper (see Remark 2.1).

**2.2. A model for deformation-dependent diffusion.** We define the first and second invariants of the tensor $\mathbf{E}_l$ as follows:

$$I_{\mathbf{E}_l} := \mathrm{tr}[\mathbf{E}_l] \tag{4a}$$

$$II_{\mathbf{E}_l} := \sqrt{2\,\mathrm{dev}[\mathbf{E}_l] \bullet \mathrm{dev}[\mathbf{E}_l]} = \sqrt{\frac{2}{3}\left(3\,\mathrm{tr}[\mathbf{E}_l^2] - (\mathrm{tr}[\mathbf{E}_l])^2\right)} \tag{4b}$$

where $\mathrm{dev}[\mathbf{E}_l] := \mathbf{E}_l - \frac{1}{3}\mathrm{tr}[\mathbf{E}_l]\mathbf{I}$ is the deviatoric part of $\mathbf{E}_l$. (Note that the second invariant is not the principal invariant, and the reason for such a choice will be discussed later.) For a given strain field (that is, for a given deformation field), we model the effect of deformation on the diffusivity as follows:

$$\mathbf{D}_{\mathbf{E}_l}(\mathbf{x}) = \mathbf{D}_0(\mathbf{x}) + (\mathbf{D}_T(\mathbf{x}) - \mathbf{D}_0(\mathbf{x}))\left(\frac{\exp[\eta_T I_{\mathbf{E}_l}] - 1}{\exp[\eta_T E_{\mathrm{ref}}] - 1}\right) + (\mathbf{D}_S(\mathbf{x}) - \mathbf{D}_0(\mathbf{x}))\left(\frac{\exp[\eta_S II_{\mathbf{E}_l}] - 1}{\exp[\eta_S E_{\mathrm{ref}}] - 1}\right) \tag{5}$$

where $\eta_T$ and $\eta_S$ are non-negative parameters; $\mathbf{D}_0(\mathbf{x})$, $\mathbf{D}_T(\mathbf{x})$ and $\mathbf{D}_S(\mathbf{x})$ are (respectively) the reference diffusivity tensors under no, tensile, and shear strains; and $E_{\mathrm{ref}}$ is a reference measure of the strain. The above model is partly motivated by the stress-induced diffusion experiments on glass done by McAfee [38, 39]. These experiments have clearly shown the following aspects, which have been qualitatively incorporated in the above model.

(a) The relative diffusion rate under tension is nearly five times more than that of the relative diffusion rates under compression and shear.
(b) The relative diffusion rate varies exponentially with respect to the (circumferential) strain for Pyrex glass (see [38, Equation 15 and Figure 3]).
(c) The relative diffusion rate under compression is significantly different from that of shear (see [39, Figure 4]).

**Remark 2.1.** *In Reference [25], the effect of deformation on the diffusivity tensor is modeled using the Frobenius norm of the Almansi-Hamel strain. (In linearized elasticity, the Almansi-Hamel strain is approximately equal to linearized strain.) Although the Frobenius norm of the strain is an invariant, the model cannot capture the difference in diffusivity under tension, compression, and shear, which has been observed in many materials. On the other hand, our model given by equation (5) can capture such departures between the diffusivities under tension, compression and shear.*

**Remark 2.2.** *A remark on the choice of invariants in the model (given by equation (5)) is in order. Note that the second principal invariant of a tensor $\mathbf{A}$ is defined as*

$$II_{\mathbf{A}}^* := \frac{1}{2}\left((\mathrm{tr}[\mathbf{A}])^2 - \mathrm{tr}[\mathbf{A}^2]\right) \tag{6}$$

*which is different from the one used in equation (4). It has been discussed in the literature that the principal invariants are not suitable to fit experimental data (e.g., Lurie [33], Anand [2, 3], Criscione et al. [13], Plešek and Kruisová [49]). These works employed invariants of Eulerian Hencky strain. Let the Eulerian Hencky strain be denoted by $\mathbf{E}_H = \ln[\mathbf{V}]$, with $\mathbf{V}$ being the left stretch tensor in the polar decomposition of the deformation gradient $\mathbf{F}$. The first three invariants based on the Eulerian Hencky strain which represent dilation ($k_1$), magnitude of distortion ($k_2$) and mode of distortion ($k_3$) are given by*

$$k_1 = \mathrm{tr}[\mathbf{E}_H] = \ln[J] \tag{7a}$$

$$k_2 = \sqrt{\mathrm{dev}[\mathbf{E}_H] \bullet \mathrm{dev}[\mathbf{E}_H]} \tag{7b}$$

$$k_3 = 3\sqrt{6}\det\left[\frac{1}{k_2}\mathrm{dev}[\mathbf{E}_H]\right] \tag{7c}$$

where $\text{dev}[\mathbf{E}_H] := \mathbf{E}_H - \frac{1}{3}\text{tr}[\mathbf{E}_H]\mathbf{I}$ and $J = \det[\mathbf{F}] > 0$. For small gradients of the displacement, we have

$$\mathbf{E}_H \approx \mathbf{E}_l \qquad (8)$$

In linearized elasticity, the constitutive equation depends only on $k_1$ and $k_2$ as any dependence on $k_3$ makes the model inherently non-linear (see Criscione et al. [13]). Herein, for modeling deformation-dependent diffusivity tensor we did not use $k_3$ just to be consistent with the theory of linearized elasticity. However, in the case of finite elasticity, one should also use the third invariant in modeling the deformation-dependent diffusivity.

One of the attractive features of the above model is that the parameters in the model can be established using standard experiments, which we shall describe below.

(a) It is easy to check that if there is no strain (that is, $\mathbf{E}_l = \mathbf{0}$), then $\mathbf{D}_{\mathbf{E}_l} = \mathbf{D}_0$. Hence, one can find the reference diffusivity tensor $\mathbf{D}_0$ by doing a diffusion experiment on an unstrained specimen.

(b) Under simple shear, we have $I_{\mathbf{E}_l} = 0$, and $II_{\mathbf{E}_l} = \alpha_S$, where $\alpha_S$ is the angle of shear. As $\alpha_S \to E_{\text{ref}}$, the diffusivity tensor $\mathbf{D} \to \mathbf{D}_S$. Hence, one can find $\mathbf{D}_S$ and $\eta_S$ by doing a diffusion experiment on a specimen undergoing a simple shear.

(c) Under a uniform tri-axial tension test, the linearized strain can be written as $\mathbf{E}_l = \alpha_T \mathbf{I}$ with $\alpha_T > 0$. In this case, $II_{\mathbf{E}_l} = 0$, $I_{\mathbf{E}_l} = 3\alpha_T$, and as $\alpha_T \to \frac{1}{3}E_{\text{ref}}$ we have $\mathbf{D} \to \mathbf{D}_T$. Hence, one can find $\mathbf{D}_T$ and $\eta_T$ in the model given by equation (5) by doing a diffusion experiment on a specimen under uniform tri-axial strain.

There are experimental techniques discussed in the literature for maintaining a specimen under uniform tri-axial strain. To name a few, tri-axial tensile testing of brittle materials such as calestone (a dental plaster), copper and aluminum alloys, austenitic stainless steel were described by Cridland and Wood [12], Hayhurst and Felce [21] and Calloch and Marquis [7]. Advanced tri-axial testing of geomaterials such as rock and soil were carried out by Donaghe et al. [14], Hunsche [23] and Wawersik [54]. Despite these experimental techniques, it can well be argued that maintaining a specimen under uniform tri-axial strain can be a difficult and expensive. In that case, after determining $\mathbf{D}_0$, $\mathbf{D}_S$ and $\eta_S$ one can evaluate $\mathbf{D}_T$ and $\eta_T$ using either uni-axial or bi-axial tension tests. However, in uni-axial and bi-axial tests, it should be noted that $II_{\mathbf{E}_l}$ will not be equal to zero.

**Remark 2.3.** *For brittle materials such as concrete, ceramics, metallic alloys and geomaterials such as soil and rock it is relatively easier to perform compression tests than tension tests. In those* cases, the expression for diffusivity tensor $\mathbf{D}_{\mathbf{E}_l}(\mathbf{x})$ can be modeled as follows:

$$\mathbf{D}_{\mathbf{E}_l}(\mathbf{x}) = \mathbf{D}_0(\mathbf{x}) + (\mathbf{D}_0(\mathbf{x}) - \mathbf{D}_C(\mathbf{x})) \left( \frac{\exp[\eta_C I_{\mathbf{E}_l}] - 1}{\exp[\eta_C E_{\text{ref}}] - 1} \right) + (\mathbf{D}_S(\mathbf{x}) - \mathbf{D}_0(\mathbf{x})) \left( \frac{\exp[\eta_S II_{\mathbf{E}_l}] - 1}{\exp[\eta_S E_{\text{ref}}] - 1} \right) \tag{9}$$

where $\eta_C$ and $\eta_S$ are non-negative parameters; $\mathbf{D}_0(\mathbf{x})$, $\mathbf{D}_C(\mathbf{x})$, and $\mathbf{D}_S(\mathbf{x})$ are (respectively) reference diffusivity tensors under no, compressive, and shear strains. Also, it should be noted that modeling diffusivity as an exponential function of stress/strain is quite popular in literature (e.g., see McAfee [38], Fahmy and Hurt [16]).

The diffusivity tensor is assumed to be symmetric, bounded, and uniformly elliptic. That is,

$$\mathbf{D}_{\mathbf{E}_l}(\mathbf{x}) = \mathbf{D}_{\mathbf{E}_l}^T(\mathbf{x}) \quad \forall\, \mathbf{x} \in \Omega \tag{10}$$

and there exists two constants $0 < \xi_1 \leq \xi_2 < +\infty$ such that

$$\xi_1 \mathbf{y}^T \mathbf{y} \leq \mathbf{y}^T \mathbf{D}_{\mathbf{E}_l}(\mathbf{x}) \mathbf{y} \leq \xi_2 \mathbf{y}^T \mathbf{y} \quad \forall\, \mathbf{x} \in \Omega \ \text{and} \ \forall\, \mathbf{y} \in \mathbb{R}^{nd} \tag{11}$$

**2.3. Governing field equations.** The governing equations for the deformation of the solid can be written as follows:

$$\text{div}[\mathbf{T}_c] + \rho(\mathbf{x})\mathbf{b}(\mathbf{x}) = \mathbf{0} \quad \text{in } \Omega \tag{12a}$$

$$\mathbf{u}(\mathbf{x}) = \mathbf{u}^{\text{P}}(\mathbf{x}) \quad \text{on } \Gamma_u^{\text{D}} \tag{12b}$$

$$\mathbf{T}_c\, \mathbf{n}(\mathbf{x}) = \mathbf{t}^{\text{P}}(\mathbf{x}) \quad \text{on } \Gamma_u^{\text{N}} \tag{12c}$$

where $\rho(\mathbf{x})$ is the density, $\mathbf{b}(\mathbf{x})$ is the specific body force, $\mathbf{u}^{\text{P}}(\mathbf{x})$ is the prescribed displacement, $\mathbf{t}^{\text{P}}(\mathbf{x})$ is the prescribed traction, and recall that $\mathbf{n}(\mathbf{x})$ is the unit outward normal to the boundary. In the absence of internal couples, the balance of angular momentum reads

$$\mathbf{T}_c = \mathbf{T}_c^T \tag{13}$$

which the Cauchy stress given in equation (2) clearly satisfies. The governing equations for the steady-state (deformation-dependent) diffusion process can be written as follows:

$$-\text{div}\left[ \mathbf{D}_{\mathbf{E}_l}(\mathbf{x}) \, \text{grad}[c] \right] = f(\mathbf{x}) \quad \text{in } \Omega \tag{14a}$$

$$c(\mathbf{x}) = c^{\text{P}}(\mathbf{x}) \quad \text{on } \Gamma_c^{\text{D}} \tag{14b}$$

$$\mathbf{n}(\mathbf{x}) \cdot \mathbf{D}_{\mathbf{E}_l}(\mathbf{x}) \, \text{grad}[c] = h^{\text{P}}(\mathbf{x}) \quad \text{on } \Gamma_c^{\text{N}} \tag{14c}$$

where $c^{\text{P}}(\mathbf{x})$ is the prescribed concentration, $h^{\text{P}}(\mathbf{x})$ is the prescribed concentration flux on the boundary, and $f(\mathbf{x})$ is the volumetric source.

It is easy to see that the governing equations for the deformation (12) and the governing equations for the diffusion (14) are coupled through equations (2) and (5). To predict the degradation/healing of the solid due to the diffusant, one needs to solve this system of coupled nonlinear partial differential equations given by equations (1)-(4), (5), (12) and (14). It is noteworthy that, except for simple problems, it is not possible to find analytical solutions to this system of equations, and one needs to resort to numerical solutions for solving realistic and practical problems.

However, to obtain reliable and predictive numerical solutions, one has to overcome many numerical challenges. In particular, one has to make sure that the chosen numerical scheme gives non-negative values for the concentration as a negative value for the concentration is unphysical. We will show in a subsequent section that the classical single-field formulation produces negative concentrations. This is particularly true if the medium is anisotropic, and one may even get negative values for the concentration even when the medium is isotropic if the mesh is not chosen with care. For example, one need to choose a mesh with square elements or a well-centered triangular mesh [45] even for isotropic diffusion. (In two-dimensions a well-centered triangular mesh means that all the angles of any triangle are acute. Similarly, one can define a well-centered mesh in higher dimensions [53].)

*We now present a numerical framework to solve the coupled deformation-diffusion equations* (1) *–*(14) *in a systematic manner.*

## 3. A COUPLED COMPUTATIONAL FRAMEWORK

The computational framework will be based on the Finite Element Method (FEM), convex quadratic programming, and staggered coupling techniques. To this end, let the domain $\Omega$ be decomposed into "$Nele$" non-overlapping open sub-domains (which in the finite element context will be elements). That is,

$$\bar{\Omega} = \bigcup_{e=1}^{Nele} \bar{\Omega}^e \quad (15)$$

where a superposed bar denotes the set closure. The boundary of $\Omega^e$ is denoted as $\partial \Omega^e := \bar{\Omega}^e - \Omega^e$. For convenience and to avoid errors due to projection operators, we shall employ the same computational mesh for both deformation and diffusion analyses. (Note that one needs to employ projection operators if different computational meshes are employed for multi-field problems like coupled deformation-diffusion.) We now present individual solvers for deformation and diffusion, and a coupling algorithm to couple individual solvers to obtain the coupled response. The solver for the diffusion problem will always give physically meaningful non-negative values for the concentration.

3.1. **A numerical solver for deformation analysis.** For solving the pure deformation problem (that is, for a given concentration field) the computational framework utilizes the standard single-field (pure displacement) formulation. However, it should be noted that one can employ any other formulation (e.g., B-bar method [22], stabilized mixed formulation [43], mixed assumed strain formulation [50], mixed enhanced formulation [52]) to solve the (pure) deformation problem. For completeness, we shall now briefly outline the single-field formulation. To this end, we shall define the following function spaces:

$$\mathcal{U} := \left\{ \mathbf{u}(\mathbf{x}) \in \left(H^1(\Omega)\right)^{nd} \mid \mathbf{u}(\mathbf{x}) = \mathbf{u}^{\mathrm{P}}(\mathbf{x}) \text{ on } \Gamma_u^{\mathrm{D}} \right\} \tag{16a}$$

$$\mathcal{W} := \left\{ \mathbf{w}(\mathbf{x}) \in \left(H^1(\Omega)\right)^{nd} \mid \mathbf{w}(\mathbf{x}) = \mathbf{0} \text{ on } \Gamma_u^{\mathrm{D}} \right\} \tag{16b}$$

where $H^1(\Omega)$ is a standard Sobolev space on $\Omega$ [6], and recall that "$nd$" is the number of spatial dimensions. The standard single-field formulation for the pure deformation problem (12) reads: Find $\mathbf{u}(\mathbf{x}) \in \mathcal{U}$ such that we have

$$\mathcal{B}_u(\mathbf{w}; \mathbf{u}) = L_u(\mathbf{w}) \quad \forall\, \mathbf{w}(\mathbf{x}) \in \mathcal{W} \tag{17}$$

where the bilinear form and linear functional are, respectively, defined as follows:

$$\mathcal{B}_u(\mathbf{w}; \mathbf{u}) := \int_\Omega \mathrm{grad}[\mathbf{w}] \bullet \mathbf{T}_c(\mathbf{u}, \mathbf{x})\, \mathrm{d}\Omega \tag{18a}$$

$$L_u(\mathbf{w}) := \int_\Omega \mathbf{w}(\mathbf{x}) \cdot \rho(\mathbf{x}) \mathbf{b}(\mathbf{x})\, \mathrm{d}\Omega + \int_{\Gamma_u^{\mathrm{N}}} \mathbf{w}(\mathbf{x}) \cdot \mathbf{t}^{\mathrm{P}}(\mathbf{x})\, \mathrm{d}\Gamma \tag{18b}$$

For a non-negative integer $m$, let $\mathbb{P}^m(\Omega^e)$ denotes the linear vector space spanned by polynomials up-to $m$th order defined on the sub-domain $\Omega^e$. We shall define the following finite dimensional subsets of $\mathcal{U}$ and $\mathcal{W}$:

$$\mathcal{U}^h := \left\{ \mathbf{u}^h(\mathbf{x}) \in \mathcal{U} \mid \mathbf{u}^h(\mathbf{x}) \in \left(C^0(\bar{\Omega})\right)^{nd},\, \mathbf{u}^h(\mathbf{x})|_{\Omega^e} \in \left(\mathbb{P}^k(\Omega^e)\right)^{nd},\, e = 1, \cdots, Nele \right\} \tag{19a}$$

$$\mathcal{W}^h := \left\{ \mathbf{w}^h(\mathbf{x}) \in \mathcal{W} \mid \mathbf{w}^h(\mathbf{x}) \in \left(C^0(\bar{\Omega})\right)^{nd},\, \mathbf{w}^h(\mathbf{x})|_{\Omega^e} \in \left(\mathbb{P}^k(\Omega^e)\right)^{nd},\, e = 1, \cdots, Nele \right\} \tag{19b}$$

where $k$ is a non-negative integer. A corresponding finite element formulation for the deformation analysis can be written as: Find $\mathbf{u}^h(\mathbf{x}) \in \mathcal{U}^h$ such that we have

$$\mathcal{B}_u(\mathbf{w}^h; \mathbf{u}^h) = L_u(\mathbf{w}^h) \quad \forall\, \mathbf{w}^h(\mathbf{x}) \in \mathcal{W}^h \tag{20}$$

After the finite element discretization, the deformation analysis will involve solving the following discrete equations:

$$\boldsymbol{K}_u(\boldsymbol{c})\boldsymbol{u} = \boldsymbol{f}_u \tag{21}$$

where $\boldsymbol{u}$ denotes the nodal displacements, $\boldsymbol{c}$ denotes the nodal concentrations, and the stiffness matrix for the deformation analysis $\boldsymbol{K}_u$ depends on the concentration of the diffusant.

3.2. **A numerical solver for diffusion analysis.** Before we provide a numerical solver for diffusion, we shall provide a mathematical argument to show that $c(\mathbf{x}) \geq 0$ in $\Omega$ even for the coupled problem. We shall assume that a solution exists for the coupled problem. (Proving existence of a solution for the coupled system is a research topic by itself, and is beyond the scope of this paper. However, in this paper we do find a numerical solution for various coupled deformation-diffusion problems.) That is, there exists a pair, $\mathbf{u}(\mathbf{x})$ and $c(\mathbf{x})$, such that they satisfy the coupled system of equations. For the given displacement field, $\mathbf{u}(\mathbf{x})$, (and hence for a given strain field $\mathbf{E}_l(\mathbf{x})$) we define

$$\mathbf{D}(\mathbf{x}) := \mathbf{D}_{\mathbf{E}_l}(\mathbf{x}) \tag{22}$$

From the theory of partial differential equations we have the following maximum principle [15]:

**Theorem 3.1.** *Let $c^{\mathrm{p}}(\mathbf{x}) \geq 0$ on $\partial \Omega$ and $\mathbf{D}(\mathbf{x})$ be continuously differentiable. If $c(\mathbf{x}) \in C^2(\Omega) \cap C^0(\bar{\Omega})$ satisfies the differential inequality $-\mathrm{div}[\mathbf{D}(\mathbf{x})\mathrm{grad}[c]] = f(\mathbf{x}) \geq 0$ in $\Omega$, then we have the following non-negative property:*

$$c(\mathbf{x}) \geq 0 \quad \text{in } \bar{\Omega} \tag{23}$$

**Remark 3.2.** *The above maximum principle theorem is due to Hopf, and a proof can be found in any standard textbook on partial differential equations (e.g., References [40, 15, 47, 18, 20]). One can find in the literature maximum principles for diffusion-type equations under weaker regularity than $C^2(\Omega) \cap C^0(\bar{\Omega})$ (for example, see References [37, 4]). But such a thorough treatment is beyond the scope of this paper, and is not crucial to our presentation.*

As discussed earlier, many existing numerical formulations (including the single-field formulation, which is based on the Galerkin principle) for diffusion-type equation do not meet the non-negative constraint. For example, the widely used single-field formulation (which is based on the Galerkin principle) does not produce physically meaningful non-negative solutions. We now present a novel methodology of enforcing the non-negative constraint on the concentration of the diffusant. The methodology works well for general computational grids and for low-order finite elements.

We start with the single-field formulation, and then modify the underlying (discrete) variational statement to meet the non-negative constraint. We shall define the following function spaces

$$\mathcal{P} := \left\{ c(\mathbf{x}) \in H^1(\Omega) \mid c(\mathbf{x}) = c^{\mathrm{p}}(\mathbf{x}) \text{ on } \Gamma^{\mathrm{D}}_c \right\} \tag{24a}$$

$$\mathcal{Q} := \left\{ q(\mathbf{x}) \in H^1(\Omega) \mid q(\mathbf{x}) = 0 \text{ on } \Gamma^{\mathrm{D}}_c \right\} \tag{24b}$$

where $H^1(\Omega)$ is a standard Sobolev space [6]. We also relax the regularity of the diffusivity tensor for weak solutions, and assume that

$$\int_\Omega \mathrm{tr}[\mathbf{D}_{\mathbf{E}_l}(\mathbf{x})^T \mathbf{D}_{\mathbf{E}_l}(\mathbf{x})] \, \mathrm{d}\Omega < +\infty \tag{25}$$

where $\mathrm{tr}[\cdot]$ is the standard trace operator used in tensor algebra and continuum mechanics [10]. The standard single field formulation for tensorial diffusion equation (14) reads: Find $c(\mathbf{x}) \in \mathcal{P}$ such that we have

$$\mathcal{B}_c(q; c) = L_c(q) \quad \forall \, q(\mathbf{x}) \in \mathcal{Q} \tag{26}$$

where the bilinear form and linear functional are, respectively, defined as

$$\mathcal{B}_c(q; c) := \int_\Omega \mathrm{grad}[q] \cdot \mathbf{D}_{\mathbf{E}_l}(\mathbf{x}) \mathrm{grad}[c] \, \mathrm{d}\Omega \tag{27a}$$

$$L_c(q) := \int_\Omega q(\mathbf{x}) \, f(\mathbf{x}) \, \mathrm{d}\Omega + \int_{\Gamma_c^\mathrm{N}} q(\mathbf{x}) \, h^\mathrm{p}(\mathbf{x}) \, \mathrm{d}\Gamma \tag{27b}$$

It is well-known that the above weak form (26) has a corresponding variational statement, which can be written as follows:

$$\underset{c(\mathbf{x}) \in \mathcal{P}}{\mathrm{minimize}} \quad \frac{1}{2}\mathcal{B}_c(c; c) - L_c(c) \tag{28}$$

We shall define the following finite dimensional vector spaces of $\mathcal{P}$ and $\mathcal{Q}$:

$$\mathcal{P}^h := \left\{ c^h(\mathbf{x}) \in \mathcal{P} \mid c^h(\mathbf{x}) \in C^0(\bar{\Omega}), \, c^h(\mathbf{x})\big|_{\Omega^e} \in \mathbb{P}^k(\Omega^e), \, e = 1, \cdots, Nele \right\} \tag{29a}$$

$$\mathcal{Q}^h := \left\{ q^h(\mathbf{x}) \in \mathcal{Q} \mid q^h(\mathbf{x}) \in C^0(\bar{\Omega}), \, q^h(\mathbf{x})\big|_{\Omega^e} \in \mathbb{P}^k(\Omega^e), \, e = 1, \cdots, Nele \right\} \tag{29b}$$

where $k$ is a non-negative integer. (Recall that, for a non-negative integer $m$, $\mathbb{P}^m(\Omega^e)$ denotes the linear vector space spanned by polynomials up-to $m$th order defined on the subdomain $\Omega^e$.) A corresponding finite element formulation can be written as: Find $c^h(\mathbf{x}) \in \mathcal{P}^h$ such that we have

$$\mathcal{B}_c(q^h; c^h) = L_c(q^h) \quad \forall \, q^h(\mathbf{x}) \in \mathcal{Q}^h \tag{30}$$

3.2.1. *A non-negative solver for tensorial diffusion equation.* We shall use the symbols $\preceq$ and $\succeq$ to denote component-wise inequalities for vectors. That is, for given any two (finite dimensional) vectors $\boldsymbol{a}$ and $\boldsymbol{b}$

$$\boldsymbol{a} \preceq \boldsymbol{b} \quad \text{means that} \quad a_i \leq b_i \, \forall \, i \tag{31}$$

Similarly, one can define the symbol $\succeq$. We shall denote the standard inner-product on Euclidean spaces by $\langle \cdot ; \cdot \rangle$.

After finite element discretization, for given nodal displacement vector $\boldsymbol{u}$, the discrete equation for the diffusion analysis takes the following form:

$$\boldsymbol{K}_c(\boldsymbol{u})\boldsymbol{c} = \boldsymbol{f}_c \tag{32}$$

where $\boldsymbol{K}_c$ is a symmetric positive definite matrix, $\boldsymbol{c}$ is the vector containing nodal concentrations, and $\boldsymbol{f}_c$ is the source vector. Let "$ndofs$" denote the number of degrees-of-freedom for the concentration. The matrix $\boldsymbol{K}_c$ is of size $ndofs \times ndofs$, and the vectors $\boldsymbol{c}$ and $\boldsymbol{f}_c$ are of size $ndofs \times 1$. The finite element discretized equation (32) is *equivalent* to the following minimization problem:

$$\underset{\boldsymbol{c} \in \mathbb{R}^{ndofs}}{\text{minimize}} \quad \frac{1}{2} \langle \boldsymbol{c}; \boldsymbol{K}_c(\boldsymbol{u})\boldsymbol{c} \rangle - \langle \boldsymbol{c}; \boldsymbol{f}_c \rangle \tag{33}$$

As shown in Figure 15(a), the numerical formulation based on equations (32) and (33) produces unphysical negative concentrations for many practically important problems. (More examples showing Galerkin formulation producing negative solutions can be found in Reference [42].) Following the ideas outlined in References [45, 32, 42] a non-negative formulation corresponding to (33) can be written as follows:

$$\underset{\boldsymbol{c} \in \mathbb{R}^{ndofs}}{\text{minimize}} \quad \frac{1}{2} \langle \boldsymbol{c}; \boldsymbol{K}_c(\boldsymbol{u})\boldsymbol{c} \rangle - \langle \boldsymbol{c}; \boldsymbol{f}_c \rangle \tag{34a}$$

$$\text{subject to} \quad \boldsymbol{c} \succeq \boldsymbol{0} \tag{34b}$$

Since for a given nodal displacement vector, the matrix $\boldsymbol{K}_c(\boldsymbol{u})$ is positive definite, the above problem (34) belongs to convex quadratic programming. From optimization theory it can be shown that the problem (34) has a unique global minimizer. There are many robust numerical algorithms available in the literature to solve the aforementioned constrained minimization problem (e.g., active set strategy, interior point methods, barrier methods). In all our numerical simulations, we have employed the active set strategy. A detailed discussion on numerical optimization can be found in references [46, 5, 34].

**3.3. A staggered coupling algorithm.** Current coupling algorithms are broadly classified into main classes: monolithic and staggered schemes. In *monolithic schemes* discretization scheme is applied to the full problem. On the other hand, in *staggered schemes* (which are based on operator-split techniques) the coupled system is partitioned, often according to the different coupled fields (in our case, concentration and displacement in deformation-diffusion analysis), and each partition is treated by different and tailored numerical schemes. There is no easy way to incorporate our non-negative formulation within a monolithic scheme. Therefore, we shall employ a staggered coupling approach. The various steps of the proposed coupling algorithm is presented in Algorithm 1.

**Algorithm 1** Staggered coupling algorithm for deformation-diffusion analysis
1: Input: tolerance ($\epsilon_{\mathrm{TOL}}^{(c)}$), maximum number of iteration (MAXITERS)
2: Guess $\boldsymbol{c}^{(0)} \succeq \boldsymbol{0}$
3: **for** $i = 1, 2, \cdots$ **do**
4:   **if** $i >$ MAXITERS **then**
5:     Solution did not converge in specified maximum number of iterations. EXIT
6:   **end if**
7:   Call deformation solver: Obtain $\boldsymbol{u}^{(i)}$ by solving
$$\boldsymbol{K}_u(\boldsymbol{c}^{(i-1)})\boldsymbol{u}^{(i)} = \boldsymbol{f}_u$$
8:   Call non-negative diffusion solver: Obtain $\boldsymbol{c}^{(i)}$ by solving the following minimization problem
$$\underset{\boldsymbol{c}^{(i)} \in \mathbb{R}^{ndofs}}{\text{minimize}} \quad \frac{1}{2}\langle \boldsymbol{c}^{(i)}; \boldsymbol{K}_c(\boldsymbol{u}^{(i)})\boldsymbol{c}^{(i)} \rangle - \langle \boldsymbol{c}^{(i)}; \boldsymbol{f}_c \rangle$$
$$\text{subject to} \quad \boldsymbol{c}^{(i)} \succeq \boldsymbol{0}$$
9:   **if** $\|\boldsymbol{c}^{(i)} - \boldsymbol{c}^{(i-1)}\| < \epsilon_{\mathrm{TOL}}^{(c)}$ **then**
10:    Staggered scheme converged. EXIT
11:   **end if**
12: **end for**

**Remark 3.3.** *In this paper, we shall use the 2-norm in the stopping criterion* $\|\boldsymbol{c}^{(i)} - \boldsymbol{c}^{(i-1)}\| < \epsilon_{\mathrm{TOL}}^{(c)}$ *in the staggered coupling algorithm. However, one could use any other norm as in (finite dimensional) Euclidean spaces all norms are* equivalent [19].

**Remark 3.4.** *It should be noted that even though the solution procedure is a staggered scheme, the problem is coupled, and the converged numerical solution will be the coupled response.*

## 4. PERFORMANCE OF THE PROPOSED FRAMEWORK

In this section, we present numerical solutions of several realistic problems using the proposed mathematical model and computational framework. The first set of problems involve degradation/healing of beams, and the second set involves degradation/healing of rectangular domains with holes. All these problems naturally arise in many important engineering applications. For example, study of degradation of beams is important to assess the performance of structural components in bridges, towers, buildings and aircraft wings. The second set of problems is to study the degradation of ducts carrying chemicals (e.g., water, carbon-dioxide, coolant) cutting through a slab or a wall, which have many applications in civil and nuclear structures. Using these two sets of representative problems we will illustrate the performance of the proposed computational framework for coupled deformation-diffusion analyses. In particular, we will show that the proposed computational framework produces physical and reliable solutions. We also systematically study the effect of the concentration of the diffusant on the deformation of the solid and vice-versa.

Here, in these representative numerical examples we have chosen $\mathbf{D}_0$, $\mathbf{D}_T$ and $\mathbf{D}_S$ as follows:

$$\mathbf{D}_0 = \begin{pmatrix} \cos(\theta) & -\sin(\theta) \\ \sin(\theta) & \cos(\theta) \end{pmatrix} \begin{pmatrix} d_1 & 0 \\ 0 & d_2 \end{pmatrix} \begin{pmatrix} \cos(\theta) & \sin(\theta) \\ -\sin(\theta) & \cos(\theta) \end{pmatrix} \tag{35a}$$

$$\mathbf{D}_T = \Phi_T \mathbf{D}_0 \tag{35b}$$

$$\mathbf{D}_S = \Phi_S \mathbf{D}_0 \tag{35c}$$

where $\Phi_T$ and $\Phi_S$ are some positive real numbers. The reason for choosing such a form is that a change in each of the $\mathbf{D}_T$ and $\mathbf{D}_S$ effect the concentration in a considerable manner in various important realistic problems. But in general each component of $\mathbf{D}_T$, $\mathbf{D}_S$ may be different from that of $\mathbf{D}_0$, even then our computational framework is still valid and works as shown in the representative problems outline in subsections plate with a hole 4.5.1 and beam with three holes 4.5.2, where in $\mathbf{D}_T$ and $\mathbf{D}_S$ are chosen as follows:

$$\mathbf{D}_T = \begin{pmatrix} \cos(\theta) & -\sin(\theta) \\ \sin(\theta) & \cos(\theta) \end{pmatrix} \begin{pmatrix} d_1^T & 0 \\ 0 & d_2^T \end{pmatrix} \begin{pmatrix} \cos(\theta) & \sin(\theta) \\ -\sin(\theta) & \cos(\theta) \end{pmatrix} \tag{36a}$$

$$\mathbf{D}_S = \begin{pmatrix} \cos(\theta) & -\sin(\theta) \\ \sin(\theta) & \cos(\theta) \end{pmatrix} \begin{pmatrix} d_1^S & 0 \\ 0 & d_2^S \end{pmatrix} \begin{pmatrix} \cos(\theta) & \sin(\theta) \\ -\sin(\theta) & \cos(\theta) \end{pmatrix} \tag{36b}$$

4.1. **Numerical $h$-convergence study.** In this subsection, we perform a numerical $h$-convergence study by employing the proposed computational framework on a coupled deformation-diffusion problem. Herein, by $h$-convergence we mean the overall convergence of the proposed framework with respect to the refinement of the computational mesh (but still maintaining the same order of interpolation within each finite element). Since the non-negative solver for diffusion works only for low-order finite element, we employ low-order finite elements in all our numerical simulations. *It should be noted that for each successful coupled analysis using the proposed computational framework, the staggered coupling algorithm should converge. For each iteration in the staggered coupling algorithm, the active-set strategy in the non-negative solver for the diffusion problem should converge.*

We employ the method of manufactured solutions [27] in this convergence study. The computational domain is taken as a bi-unit square ($0 \leq x \leq 1$ and $0 \leq y \leq 1$). The displacement vector is

given by

$$u(x,y) = \frac{1}{\pi}\sin(\frac{\pi x}{2})\sin(\frac{\pi y}{2}) \tag{37a}$$

$$v(x,y) = \frac{1}{\pi}\cos(\frac{\pi x}{2})\cos(\frac{\pi y}{2}) \tag{37b}$$

The concentration is given by

$$c(x,y) = 1 + \frac{1}{\pi}\sin(\frac{\pi x}{2})\sin(\frac{\pi y}{2}) \tag{38}$$

Note that the concentration given in the above equation is non-negative in the whole computational domain. The following parameters are assumed in the convergence study:

$$\mu_0 = \pi + 2, \ \mu_1 = -\pi, \ \lambda_0 = 2, \ \lambda_1 = -1, \ \rho = 1, \ E_{\text{ref}} = 0.0001,$$

$$c_{\text{ref}} = 1, \ \eta_T = 1, \ \eta_S = 1, \ \beta_T = 2, \ \beta_S = 2, \ \gamma_S = \frac{\beta_S - 1}{\exp[\eta_S E_{\text{ref}}] - 1},$$

$$\mathbf{D}_0 = 2\,\mathbf{I}, \ \mathbf{D}_T = \beta_T \mathbf{D}_0, \ \mathbf{D}_S = \beta_S \mathbf{D}_0 \tag{39}$$

The volumetric source for the diffusion problem is given by

$$f(x,y) = \pi \sin(\frac{\pi x}{2})\sin(\frac{\pi y}{2})\left((1-\gamma_S) + \gamma_S \exp\left[\cos(\frac{\pi x}{2})\sin(\frac{\pi y}{2})\right]\right)$$
$$- \frac{\pi \gamma_S}{4}\sin(\pi x)\cos(\pi y)\exp\left[\cos(\frac{\pi x}{2})\sin(\frac{\pi y}{2})\right] \tag{40}$$

The specific body force for the deformation problem is given by

$$\mathbf{b}(x,y) = \begin{pmatrix} \pi \sin(\frac{\pi x}{2})\sin(\frac{\pi y}{2}) + \frac{\pi}{4}\cos(\pi x)(1-\cos(\pi y)) \\ \pi \cos(\frac{\pi x}{2})\cos(\frac{\pi y}{2}) - \frac{\pi}{4}\sin(\pi x)\sin(\pi y) \end{pmatrix} \tag{41}$$

The boundary value problem is pictorially described in Figure 2. The tolerance in the stopping criterion for the staggered coupling algorithm is taken as $\epsilon_{\text{TOL}}^{(c)} = 10^{-8}$. A hierarchy of meshes similar to the ones shown in Figure 3 is employed in the numerical convergence study. In Figure 4, the obtained numerical solution using the mesh shown in Figure 3(a) is compared with the analytical solution. The convergence of the proposed coupled framework is illustrated in Figures 5 and 6 for three-node triangular and four-node quadrilateral elements with respect to $L^2$-norm and $H^1$-seminorm, and the proposed computational framework performed well.

In the subsequent subsections, we employ the proposed computational framework to solve various finite domain practical problems. Using these we illustrate how the tension, compression and shear in the solid affect the (steady-state) diffusion process, and how the concentration of the diffusant affect the deformation of the solid. The values of the parameters that are common to all the test problems are presented in Table 1. We have given the values for the volumetric source and body force for all the test problems in Table 2.

TABLE 1. Values of the parameters that are the same in all test problems (except for the one in the numerical $h$-convergence).

| Parameter | value |
|---|---|
| $\rho$ | 1 |
| $c_{\text{ref}}$ | 1 |
| $E_{\text{ref}}$ | 0.0001 |
| $\lambda_0$ | $+10^6$ |
| $\mu_0$ | $+10^6$ |
| $\lambda_1$ | $-9 \times 10^5$ |
| $\mu_1$ | $-9 \times 10^5$ |

TABLE 2. Volumetric source $f(\mathbf{x})$ and body force $\mathbf{b}(\mathbf{x})$ for various test problems.

| Test problem | $f(\mathbf{x})$ | $\mathbf{b}(\mathbf{x})$ |
|---|---|---|
| Cantilever beam with edge shear | 10000 | $\mathbf{0}$ |
| Simply supported beam under self-weight | 1000 | $-10\,\hat{\mathbf{e}}_y$ |
| Fixed beam under self-weight | 100 | $-10\,\hat{\mathbf{e}}_y$ |
| Plate with a hole under self-weight | 0 | $-10\,\hat{\mathbf{e}}_y$ |
| Beam with three holes under self-weight | 0 | $-10\,\hat{\mathbf{e}}_y$ |

4.2. **Cantilever beam with edge shear.** The purpose of this test problem is to illustrate the affect of $\mathbf{D}_S$ on the coupled response. We consider a cantilever beam with a uniform edge shear of 500. The length of the beam is 1.0, and the height of the beam is 0.1. For the deformation problem, top and bottom surfaces are subjected to zero traction. For the diffusion problem, zero concentration is prescribed on the top and bottom surfaces, and zero flux is prescribed on the left and right surfaces. The boundary value problem is pictorially described in Figure 7.

We have employed a structured mesh using four-node quadrilateral elements with 21 nodes along each side of the domain (see the mesh in Figure 8). The stopping criterion is again taken as $\epsilon_{\text{TOL}}^{(c)} = 10^{-8}$. We have considered three different values for $\Phi_S$: 5, 10, and 20. The other parameters used in the coupled analysis are as follows:

$$\theta = \frac{\pi}{6},\ d_1 = 1,\ d_2 = 1,\ \Phi_T = 2,\ \eta_T = 100,\ \eta_S = 1 \tag{42}$$

The contour profiles for the concentration for these three values are shown in Figure 8. As expected, the diffusant accumulated near the right side of the domain (where the uniform edge shear is applied). In Table 3, we have presented the maximum concentration and the number of iterations taken by the staggered coupling algorithm for the chosen values of $\Phi_S$. One can observe that as $\Phi_S$ increases, the staggered coupling algorithm takes fewer iterations. This can be explained as $\Phi_S$ increases, the diffusivity increases and the diffusant is more uniformly distributed (which is evident from Figure 8).

TABLE 3. Cantilever beam with edge shear: The table presents the maximum concentration and the number of iterations taken by the staggered coupling algorithm for various values of $\Phi_S$. (Note that the minimum concentration in all the cases is zero.)

| $\Phi_S$ | Max. concentration | Iterations |
|---|---|---|
| 5 | $4.257 \times 10^{-1}$ | 14 |
| 10 | $2.187 \times 10^{-1}$ | 9 |
| 20 | $1.107 \times 10^{-1}$ | 7 |

TABLE 4. Simply supported beam under self-weight: The table presents the maximum concentration and the iterations taken by the staggered coupling algorithm for various values of $\eta_S$. (Note that the minimum concentration in all the cases is zero.)

| $\eta_S$ | Max. concentration | Iterations |
|---|---|---|
| 1 | $7.205 \times 10^{-1}$ | 10 |
| $10^3$ | $7.309 \times 10^{-1}$ | 10 |
| $2 \times 10^4$ | $9.365 \times 10^{-1}$ | 12 |

4.3. **Simply supported beam under self-weight.** The purpose of this test problem is to study the affect of $\eta_S$ on the coupled response, and to illustrate the competitive effects of shear and tension/compression on the diffusion process. We consider a simply supported beam under self-weight with traction-free surfaces. The boundary value problem is pictorially described in Figure 9.

The tolerance for the stopping criterion in the staggered coupling algorithm is taken as $\epsilon_{\text{TOL}}^{(c)} = 10^{-5}$. We have employed a structured mesh using four-node quadrilateral element with 21 nodes along each side of the domain (see the mesh in Figure 10). We have considered three values for $\eta_S$: 1, $10^3$ and $2 \times 10^4$. The other parameters used in this test problem are as follows:

$$\theta = \frac{\pi}{6}, \; d_1 = 1, \; d_2 = 1, \; \Phi_T = 10, \; \Phi_S = 10, \; \eta_T = 1 \tag{43}$$

From Table 4 one can infer that as $\eta_S$ increases the maximum concentration and iterations taken for convergence of staggered coupling algorithm increase. This is because for lower values of $\eta_S$ shear affects the diffusivity tensor more than the tension/compression and for higher values of $\eta_S$ tension/compression dominates. This is also evident in the numerical results presented in Figure

TABLE 5. Fixed beam under self-weight: The table presents the maximum concentration and the iterations taken by the staggered coupling algorithm for various values of $\eta_T$. (Note that the minimum concentration in all the cases is zero.)

| $\Phi_T$ | Max. concentration | Iterations |
|---|---|---|
| 1 | $1.250 \times 10^{-1}$ | 2 |
| 5 | $1.348 \times 10^{-1}$ | 5 |
| 7 | $1.575 \times 10^{-1}$ | 8 |

10. For lower values of $\eta_S$, the diffusant accumulated near the supports of the beam (where the shear is maximum). For higher values of $\eta_S$, the diffusant spreads deep into beam. An important feature to be noted is that the concentration profile curves up. This is because the diffusivity is higher in tension than in compression. This is physically meaningful, as, in general, the sizes of pores in the solid enlarge due to tension, shrink due to compression, and distort due to shear. The deformation-dependent diffusivity tensor $\mathbf{D}_{\mathbf{E}_l}(\mathbf{x})$ given by equation (5) takes these factors into account. The next test problem highlights other important features of the proposed mathematical model.

4.4. **Fixed beam under self-weight.** The purpose of this test problem is to study the affect of $\Phi_T$ on the coupled response and show that our model can capture the qualitative aspects of the experiments done by McAfee [38, 39]. A fixed beam of length unity and depth 0.1 is subjected to self-weight. The top and bottom surfaces are traction-free. For the diffusion problem, the top and bottom surfaces have zero concentration, and the left and right surfaces have zero flux. The boundary value problem is pictorially described in Figure 11.

The tolerance for the stopping criterion in the staggered coupling algorithm is taken as $\epsilon_{\text{TOL}}^{(c)} = 10^{-7}$. We have employed a structured mesh using four-node quadrilateral element with 21 nodes along each side of the domain (see the mesh in Figure 12). We have considered three different values for $\Phi_T$: 1, 5 and 7. The other parameters used in this test problem are as follows:

$$\theta = 0, \ d_1 = 1, \ d_2 = 1, \ \Phi_S = 1, \ \eta_T = 100, \ \eta_S = 1 \tag{44}$$

In Table 5, we have presented the maximum concentration and number of iterations taken by the staggered coupling algorithm for chosen values of $\Phi_T$. One can observe that as $\Phi_T$ increases the maximum concentration and iterations taken increase. This results in the deformation-dependent diffusivity $\mathbf{D}_{\mathbf{E}_l}$ being high in tensile and low in compressive and shear regions of the beam. The diffusant accumulates in low diffusive areas such as the supports which are under shear and top section of the beam which is in compression. From Figure 12 one can observe that the concentration profile which is initially around the central line of the beam curves up (the region of interest being around 0.12). This test problem shows the success of our mathematical model in predicting the underlying physical phenomena observed from the stress-induced experiments done by McAfee [38, 39]. The invariants which represent dilation and distortion in our model were instrumental in capturing the change in diffusivity under tension, compression and shear. Hence concentration curves towards the compression zone.

This type of curving of the concentration profile is not observed when one uses the model in which $\mathbf{D}_{\mathbf{E}_l}(\mathbf{x})$ depends only on $\|\mathbf{E}_l\|$. As mentioned in the Remark 2.1, a diffusivity model based on the frobenius norm of $\mathbf{E}_l$ cannot capture the change in diffusivity tensor under tension, compression and shear deformations. The diffusivity tensor as per Reference [25] and the parameters assumed for the numerical study are given by:

$$\mathbf{D}_{\mathbf{E}_l}(\mathbf{x}) = \mathbf{D}_0(\mathbf{x}) + (\mathbf{D}_\infty(\mathbf{x}) - \mathbf{D}_0(\mathbf{x})) \left(1 - \exp[-\lambda\|\mathbf{E}_l\|]\right) \tag{45a}$$

$$\theta = 0, \ d_1 = 1, \ d_2 = 1, \ \lambda = 100, \ \mathbf{D}_\infty = 10\,\mathbf{D}_0 \tag{45b}$$

The concentration profile as illustrated in Figure 13 is always around the central line and does not curve up for any value of $\mathbf{D}_\infty$ and $\lambda$.

4.5. **Degradation/healing of rectangular domains with holes.** In the following subsections, we highlight the importance of non-negative solutions and its impact on coupled deformation-diffusion analyses of rectangular domains with holes. We compare the non-negative formulation with the standard Galerkin formulation using two representative problems. We shall model $\mathbf{D}_T$ and $\mathbf{D}_S$ by the equations (36) in which $d_j^i$ is different from $d_j$ where $i = T, S$ and $j = 1, 2$. In both numerical studies here, the tolerance for the stopping criterion in the staggered coupling algorithm is taken as $\epsilon_{\mathrm{TOL}}^{(c)} = 10^{-5}$ and three-node triangular unstructured meshes were used (see Figures 14 and 18).

4.5.1. *Plate with a hole under self-weight.* The purpose of this test problem is to study the importance of non-negative constraint in coupled-deformation-diffusion analyses. Herein, we perform coupled analysis for a square plate with a square hole under self-weight. The computational domain $\Omega$ is the region in-between a bi-unit square plate and a square hole of length $\frac{1}{9}$. The boundary conditions for displacements at the hole and traction for the plate are equal to zero. The concentration at the hole is maintained at 1 while that at boundary of the plate is equal to 0. The boundary value problem is pictorially described in Figure 14. The other parameters assumed in the analysis of the coupled problem are as follows:

$$\theta = \frac{\pi}{3}, \ \eta_T = 1, \ \eta_S = 1, \ d_1 = 10000, \ d_1^T = 11000, \ d_1^S = 11000, \ d_2 = 1, \ d_2^T = 10, \ d_2^S = 5 \tag{46}$$

From Figure 15 it is evident that the proposed non-negative formulation satisfies the above condition and produces physically meaningful concentration while the standard Galerkin formulation predicts otherwise. The comparison of the concentration and degradation profiles for the standard Galerkin vs. the non-negative formulation given in Figures 15 and 16 shows that the Galerkin formulation violates the non-negative constraint. The obtained negative values for the concentration are not just numerical noise (see Table 6).

Table 6. Plate with a hole under self-weight: The table presents the minimum concentration, degradation index and the iterations taken by staggered coupling algorithm using the standard single-field formulation. Analysis is done for various values of $(d_1, d_1^T, d_1^S)$ by keeping other parameters in equations 46 fixed.

| $(d_1, d_1^T, d_1^S)$ | Min. concentration | Degradation Index (% of nodes violated) | Iterations (Galerkin) | Iterations (non-negative) |
|---|---|---|---|---|
| $(1, 10, 10)$ | 0 | 0 | 9 | 9 |
| $(10, 30, 20)$ | 0 | 0 | 8 | 8 |
| $(100, 120, 110)$ | $-3.301 \times 10^{-3}$ | 25.15 | 6 | 21 |
| $(1000, 1200, 1100)$ | $-3.586 \times 10^{-2}$ | 32.76 | 5 | 22 |
| $(10000, 11000, 11000)$ | $-4.398 \times 10^{-2}$ | 34.07 | 4 | 8 |

4.5.2. *Beam with three holes under self-weight.* The purpose of this test problem is to show that as number of holes increase, the spatial extent and magnitude of violation of the non-negative constraint by the Galerkin formulation increase dramatically. Herein, we perform coupled deformation-diffusion analysis for a beam with three square holes under self-weight. The computational domain $\Omega$ is the region in-between a beam of length 10.0 and height 1.0 and three square holes each of length 0.4. The boundary conditions for displacements at the holes and traction for the beam are equal to zero. The concentration at the holes is maintained at 1 while that at boundary of the beam is equal to 0. The boundary value problem is pictorially described in Figure 17. The other parameters assumed in the analysis of the coupled problem are as follows:

$$\theta = \frac{\pi}{4}, \ \eta_T = 1, \ \eta_S = 1, \ d_1 = 10000, \ d_1^T = 20000, \ d_1^S = 15000, \ d_2 = 1, \ d_2^T = 5, \ d_2^S = 2 \quad (47)$$

From Figures 18 and 19 it is evident that as the number of holes increase, the regions of negative concentration obtained from the standard Galerkin formulation also increases. Also from Tables 6 and 7 it is evident that the negative concentration obtained from numerical study on beam with three holes is 1.35 times that of the plate with a hole.

From Tables 6 and 7, one can see that as the directional diffusivities $d_1$, $d_1^T$ and $d_1^S$ increase the minimum concentration (which is initially zero) becomes negative under the Galerkin formulation. This negative value of concentration and the degradation index (which represents the % of nodes at which concentration is negative) increases as these diffusivities increase. The proposed non-negative formulation does not give negative values for the concentration for any values for these how these diffusivities. This is quite important because the material properties of the solid are dependent on the concentration obtained. From Figures 16 and 19, one can observe that the standard Galerkin method predicts that some regions are *healing* due to the *negative* values of the concentration which is physically unrealistic. The non-negative solver employed always gives non-negative concentration and shows that material is *degrading everywhere*.

TABLE 7. Beam with three holes under self-weight: The table presents the minimum concentration, degradation index and the iterations taken by staggered coupling algorithm using the standard single-field formulation. Analysis is done for various values of $(d_1, d_1^T, d_1^S)$ by keeping other parameters in equations 47 fixed.

| $(d_1, d_1^T, d_1^S)$ | Min. concentration | Degradation Index (% of nodes violated) | Iterations (Galerkin) | Iterations (non-negative) |
| --- | --- | --- | --- | --- |
| $(1, 10, 10)$ | 0 | 0 | 8 | 8 |
| $(10, 50, 30)$ | 0 | 0 | 8 | 8 |
| $(100, 200, 150)$ | $-3.150 \times 10^{-2}$ | 30.69 | 6 | 8 |
| $(1000, 2000, 1500)$ | $-5.665 \times 10^{-2}$ | 31.14 | 6 | 6 |
| $(10000, 20000, 15000)$ | $-5.948 \times 10^{-2}$ | 31.14 | 6 | 16 |

4.6. **Performance of the staggered scheme and the active-set strategy.** The convergence histories of the staggered coupling algorithm for the aforementioned test problems are shown in Figure 20. Note that the convergence of the staggered coupling algorithm is monotonic for the test problems on degradation/healing of beams. But the convergence is not monotonic for the test problem on degradation/healing of rectangular domain with holes, which is illustrated in Figure 20(b)). Figure 21 shows the number of active-strategy iterations required for solving the constrained optimization problem for each iteration in the staggered coupling algorithm.

## 5. CONCLUSIONS

Many technologically and biologically important processes are coupled deformation-diffusion problems. Lately, there is a surge in research activity in studying coupled deformation-diffusion problems. However, in all these research efforts it has not been recognized that many popular numerical formulations and existing computational packages predict unphysical negative solutions for the concentration of the diffusant, and this is more prominent in the case of a medium which has high directional diffusivities. Concentration of a chemical is a non-negative quantity, and this property has to be preserved to obtain physically meaningful numerical solution for a coupled deformation-diffusion problem.

We proposed a mathematical model for coupled deformation-diffusion analysis, which is suitable to model degradation/healing of materials and structures. The model is *fully* coupled in the sense that the deformation process is affected by the diffusion process, and the diffusion process is in turn affected by the deformation of solid. One of the main contributions is that we have presented a robust computational framework for solving coupled deformation-diffusion problems. The computational framework includes a non-negative formulation for (tensorial) diffusion equation, a numerical solver for the deformation of the solid, and a staggered coupling algorithm. An important aspect of the computational framework is that it always produces physically meaningful non-negative values for the concentration of the diffusant on any computational grid (with low-order finite elements) even for a medium which has high directional diffusivities. We have illustrated the robustness of the computational framework on representative numerical examples. We have studied systematically the affect of the deformation on the diffusion process on various practically important problems.

## ACKNOWLEDGMENTS

The research reported herein was supported by Texas Engineering Experiment Station (TEES). This support is gratefully acknowledged. The opinions expressed in this paper are those of the authors and do not necessarily reflect that of the sponsor.

## References

[1] *ABAQUS/CAE/Standard, Version 6.8-3*. Simulia, Providence, Rhode Island, www.simulia.com, 2009.

[2] L. Anand. On Hencky's approximate strain-energy function for moderate deformations. *Journal of Applied Mechanics*, 46:78–82, 1979.

[3] L. Anand. Moderate deformations in extension-torsion of incompressible isotropic elastic materials. *Journal of the Mechanics and Physics of Solids*, 34:293–304, 1986.

[4] M. Borsuk and V. Kondratiev. *Elliptic Boundary Value Problems of Second Order in Piecewise Smooth Domains*. Elsevier Science, San Diego, USA, 2006.

[5] S. Boyd and L. Vandenberghe. *Convex Optimization*. Cambridge University Press, Cambridge, UK, 2004.

[6] F. Brezzi and M. Fortin. *Mixed and Hybrid Finite Element Methods, volume 15 of Springer series in computational mathematics*. Springer-Verlag, New York, USA, 1991.

[7] S. Calloch and D. Marquis. Triaxial tension-compression tests for multiaxial cyclic plasticity. *International Journal of Plasticity*, 15:521–549, 1999.

[8] J. L. Chaboche. Continuum damage mechanics: Part I - General Concepts. *Journal of Applied Mechanics*, 55:59–64, 1988.

[9] J. L. Chaboche. Continuum damage mechanics: Part II - Damage growth, crack initiation and crack growth. *Journal of Applied Mechanics*, 55:65–72, 1988.

[10] P. Chadwick. *Continuum Mechanics: Concise Theory and Problems*. Dover Publications, Inc., Minealo, New York, 1999.

[11] D. X. Cheng, D. N. Little, R. L. Lytton, and J. C. Holste. Moisture damage evaluation of asphalt mixtures by considering both moisture diffusion and repeated-load conditions. *Transportation Research Record: Journal of the Transportation Research Board*, (1832):42–49, 2003.

[12] L. Cridland and W. G. Wood. A hydrostatic tension test of a brittle material. *International Journal of Fracture Mechanics*, 4:277–285, 1968.

[13] J. C. Criscione, J. D. Humphrey, A. S. Douglas, and W. C. Hunter. An invariant basis for natural strain which yields orthogonal stress response terms in isotropic hyperelasticity. *Journal of the Mechanics and Physics of Solids*, 48:2445–2465, 2000.

[14] R. T. Donaghe, R. C. Chaney, and S. Marshall. *Advanced Triaxial Testing of Soil and Rock*. ASTM STP: 977, 1986.

[15] L. C. Evans. *Partial Differential Equations*. American Mathematical Society, Providence, Rhode Island, USA, 1998.

[16] A. A. Fahmy and J. C. Hurt. Stress dependence of water diffusion in epoxy resin. *Polymer Composites*, 1:77–80, 1980.

[17] L. E. Fraenkel. *An Introduction to Maximum Principles and Symmetry in Elliptic Problems*. Cambridge University Press, Cambridge, UK, 2000.

[18] D. Gilbarg and N. S. Trudinger. *Elliptic Partial Differential Equations of Second Order*. Springer, New York, USA, 2001.

[19] P. R. Halmos. *Finite-Dimensional Vector Spaces*. Springer-Verlag, New York, USA, 1993.

[20] Q. Han and F. Lin. *Elliptic Partial Differential Equations*. American Mathematical Society, Providence, Rhode Island, USA, 2000.

[21] D. R. Hayhurst and I. D. Felce. Creep rupture under tri-axial tension. *Engineering Fracture Mechanics*, 25:645–664, 1986.

[22] T. J. R. Hughes. Generalization of selective integration procedures to anisotropic and nonlinear media. *International Journal for Numerical Methods in Engineering*, 15:1413–1418, 1980.

[23] U. Hunsche. *Uniaxial and triaxial creep and failure tests on rock: Experimental technique and interpretation*, pages 1–54. Visco-Plastic Behavior of Geomaterials. Springer-Verlag, 1994.

[24] L. M. Kachanov. *Introduction to Continuum Damage Mechanics*. Kluwer Academic Publishers, Dordrecht, Netherlands, 1986.

[25] S. Karra and K. R. Rajagopal. Degradation and healing in a generalized neo-Hookean solid due to infusion of a fluid. *arXiv:1007.1038*, 2010.

[26] Y. R. Kim, D. N. Little, and R. L. Lytton. Effect of moisture damage on material properties and fatigue resistance of asphalt mixtures. *Transportation Research Record: Journal of the Transportation Research Board*, (1891):48–54, 2004.

[27] P. Knupp and K. Salari. *Verification of Computer Codes in Computational Science and Engineering*. Chapman and Hall, Boca Raton, Florida, 2003.

[28] N. Kringos, T. Scarpas, A. Copeland, and J. Youtcheff. Modelling of combined physical-mechanical moisture-induced damage in asphaltic mixes, Part 2: moisture susceptibility parameters. *International Journal of Pavement Engineering*, 9:129–151, 2008.

[29] N. Kringos, T. Scarpas, C. Kasbergen, and P. Selvadurai. Modelling of combined physical-mechanical moisture-induced damage in asphaltic mixes, Part 1: governing processes and formulations. *International Journal of Pavement Engineering*, 9:115–128, 2008.

[30] D. Kuhl, F. Bangert, and G. Meschke. Coupled chemo-mechanical deterioration of cementitious materials. Part I: Modeling. *International Journal of Solids and Structures*, 41:15–40, 2004.

[31] D. Kuhl, F. Bangert, and G. Meschke. Coupled chemo-mechanical deterioration of cementitious materials Part II: Numerical methods and simulations. *International Journal of Solids and Structures*, 41:41–67, 2004.

[32] R. Liska and M. Shashkov. Enforcing the discrete maximum principle for linear finite element solutions for elliptic problems. *Communications in Computational Physics*, 3:852–877, 2008.

[33] A. I. Lurie. *Nonlinear Theory of Elasticity*. North Holland Series in Applied Mathematics and Mechanics, Elsevier Science, Netherlands, 1990.

[34] O. L. Mangasarian. *Nonlinear Programming*. SIAM, New York, USA, 1994.

[35] S. I. Marras, I. A. Ihtiaris, N. K. Hatzitrifon, K. Sikalidis, and E. C. Aifantis. A preliminary study of stress-assisted fluid penetration in ceramic bricks. *Journal of the European Ceramic Society*, 20:243–269, 2000.

[36] A. Masud and T. J. R. Hughes. A stabilized mixed finite element method for Darcy flow. *Computer Methods in Applied Mechanics and Engineering*, 191:4341–4370, 2002.

[37] V. G. Mazja and B. A. Plamenevskii. $l_p$ estimates of solutions of elliptic boundary value problems in domains with edges. *Transactions of the Moscow Mathematical Society*, 1:49–97, 1980.

[38] K. B. McAfee. Stress-enhanced diffusion in glass I. Glass under tension and compression. *Journal of Chemical Physics*, 28:218–226, 1958.

[39] K. B. McAfee. Stress-enhanced diffusion in glass II. Glass under shear. *Journal of Chemical Physics*, 28:226–229, 1958.

[40] R. McOwen. *Partial Differential Equations: Methods and Applications*. Prentice Hall, New Jersey, USA, 1996.

[41] A. Muliana, K. R. Rajagopal, and S. C. Subramanian. Degradation of an elastic composite cylinder due to the diffusion of a fluid. *Journal of Composite Materials*, 43:1225–1249, 2009.

[42] H. Nagarajan and K. B. Nakshatrala. Enforcing the non-negativity constraint and maximum principles for diffusion with decay on general computational grids. *International Journal for Numerical Methods in Fluids*, In print, DOI: 10.1002/fld.2389, 2010.

[43] K. B. Nakshatrala, A. Masud, and K. D. Hjelmstad. On finite element formulations for nearly incompressible linear elasticity. *Computational Mechanics*, 41:547–561, 2008.

[44] K. B. Nakshatrala, D. Z. Turner, K. D. Hjelmstad, and A. Masud. A stabilized mixed finite element formulation for Darcy flow based on a multiscale decomposition of the solution. *Computer Methods in Applied Mechanics and Engineering*, 195:4036–4049, 2006.

[45] K. B. Nakshatrala and A. J. Valocchi. Non-negative mixed finite element formulations for a tensorial diffusion equation. *Journal of Computational Physics*, 228:6726–6752, 2009.

[46] J. Nocedal and S. J. Wright. *Numerical Optimization*. Springer Verlag, New York, USA, 1999.

[47] M. H. Protter and H. F. Weinberger. *Maximum Principles in Differential Equations*. Springer-Verlag, New York, USA, 1999.

[48] P. A. Raviart and J. M. Thomas. A mixed finite element method for 2nd order elliptic problems. In *Mathematical Aspects of the Finite Element Method*, pages 292–315, Springer-Verlag, New York, 1977.

[49] J. Plešek and A. Kruisová. Formulation, validation and numerical procedures for Hencky's elasticity model. *Computer and Structures*, 84:1141–1150, 2006.

[50] J. C. Simo and M. S. Rifai. A class of mixed assumed strain methods and the method of incompatible modes. *International Journal for Numerical Methods in Engineering*, 29:1595–1638, 1990.

[51] P. Sofronis. The influence of mobility of dissolved hydrogen on the elastic response of a metal. *Journal of Mechanics and Physics of Solids*, 43:1385–1407, 1995.

[52] R. L. Taylor. A mixed-enhanced formulation for tetrahedral finite elements. *International Journal for Numerical Methods in Engineering*, 47:205–227, 2000.

[53] E. Vanderzee, A. N. Hirani, D. Guoy, and E. Ramos. Well-centered triangulation. Technical Report UIUCDCS-R-2008-2936, Department of Computer Science, University of Illinois at Urbana-Champaign, February 2008. Also available as a preprint at arXiv as arXiv:0802.2108v1 [cs.CG].

[54] W. R. Wawersik, L. W. Carlson, D. J. Holcomb, and R. J. Williams. New method for true-triaxial rock testing. *International Journal of Rock Mechanics and Mining Sciences*, 34:330.e1–330.e14, 1997.

[55] Y. Weitsman. Coupled damage and moisture-transport in fiber-reinforced polymeric composites. *International Journal of Solids and Structures*, 23:1003–1025, 1987.

[56] Y. J. Weitsman. Anomalous fluid sorption in polymeric composites and its relation to fluid-induced damage. *Composities Part A: applied science and manufacturing*, 37:617–623, 2006.

[57] J. Woodtli and R. Kieselbach. Damage due to hydrogen embrittlement and stress corrosion cracking. *Engineering Failure Analysis*, 7:427–450, 2000.

[58] J. C. Wu and N. A. Peppas. Numerical simulation of anomalous penetrant diffusion in polymers. *Journal of Applied Polymer Science*, 49:1845–1856, 1993.

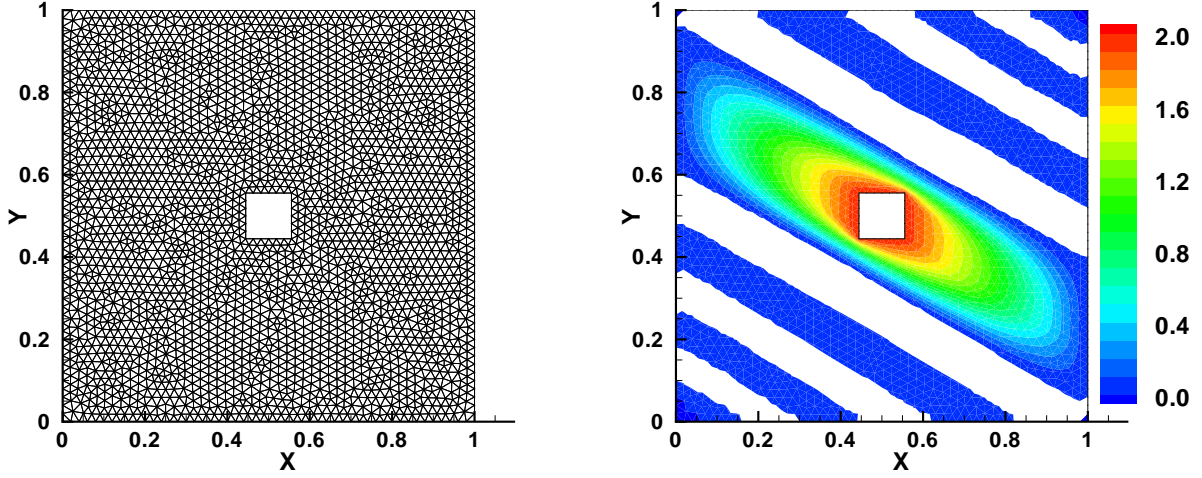

FIGURE 1. ABAQUS simulation for plate with a square hole: The left figure shows three-node triangular mesh used in the simulation using Abaqus [1]. The right figure shows contours of concentration obtained using Abaqus for the problem described in subsection 4.5.1. The white area in the right figure indicates the region in which concentration has a *negative* value. Approximately 37 % of nodes have negative values. The minimum concentration obtained is $-0.0832$ (which is -4.16 % off the range of possible values: 0 to 2). We have taken $d_1 = 10000$, $d_2 = 1$, and $\theta = -\pi/6$.

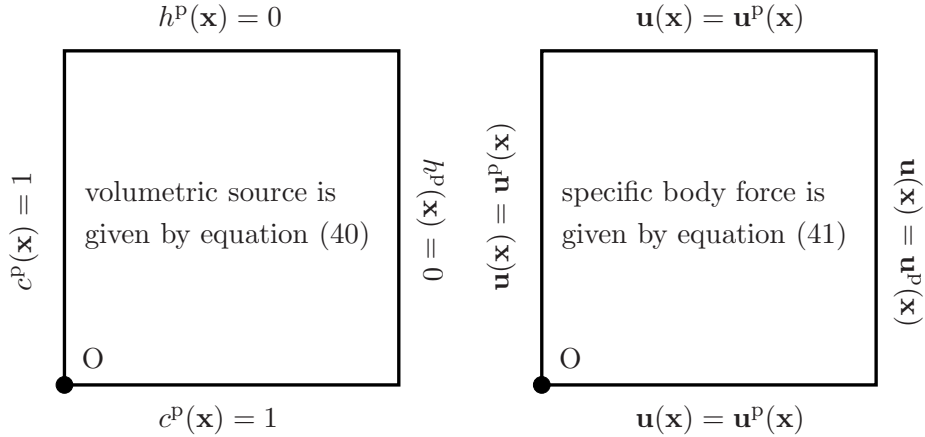

FIGURE 2. Numerical $h$-convergence: A pictorial description of the boundary value problem. The computational domain is a bi-unit square with origin denoted as O. The boundary conditions and the volumetric source for the diffusion problem is shown in the left figure. The (Dirichlet) boundary conditions and the specific body force for the deformation problem is shown in the right figure. The Dirichlet boundary conditions $\mathbf{u}^{\mathrm{p}}(\mathbf{x})$ are prescribed by directly evaluating the expressions for the displacement given in equation (37).

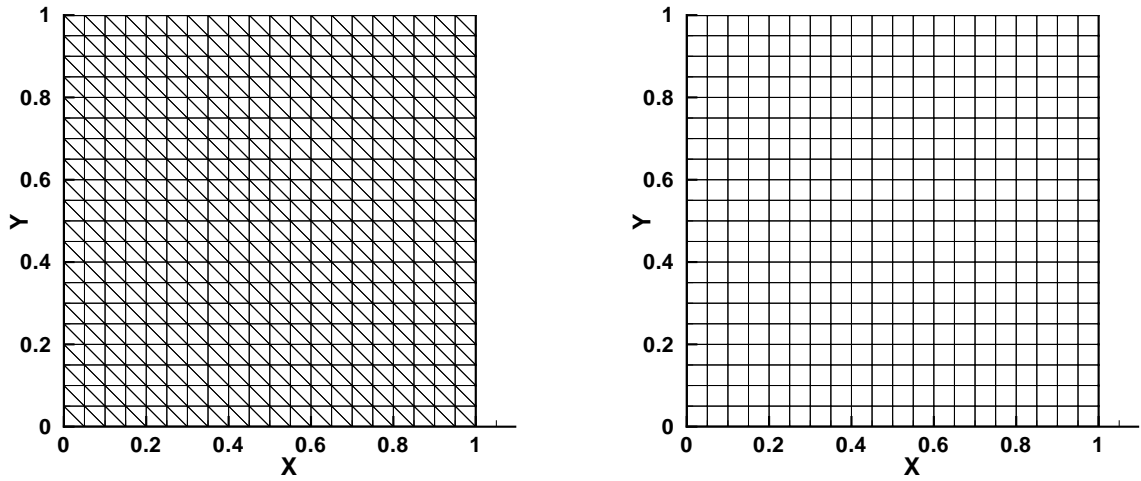

FIGURE 3. Numerical $h$-convergence: This figure shows the typical three-node triangular and four-node quadrilateral meshes used in the numerical convergence analysis. Both these meshes have 21 nodes along each side. For the numerical convergence analysis, these meshes are subdivided in a hierarchical manner to obtain a series of meshes.

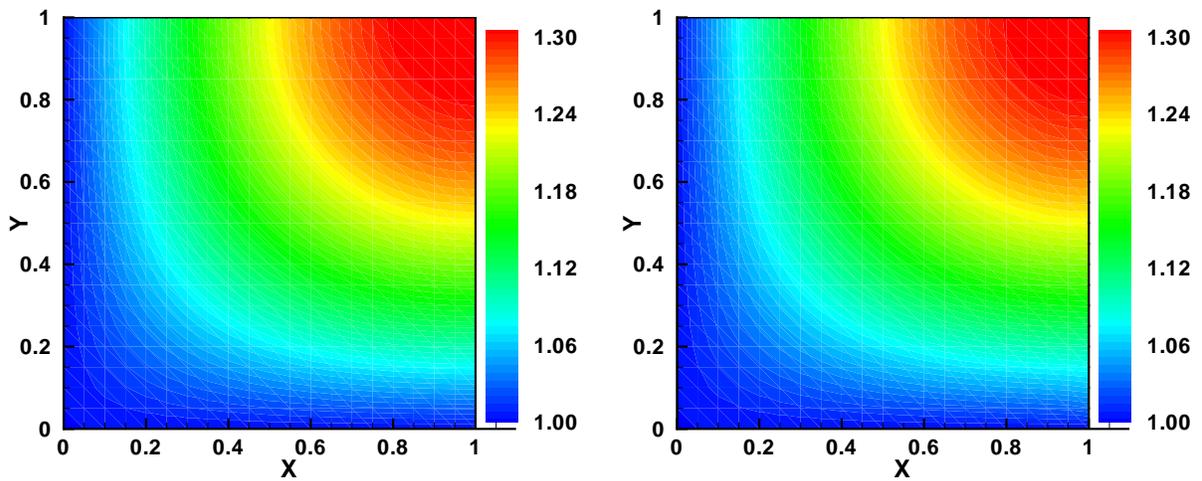

FIGURE 4. Numerical $h$-convergence: Comparison of concentration profile from analytical solution (left) to that of the numerical study (right) using three-node triangular mesh shown in Figure 3.

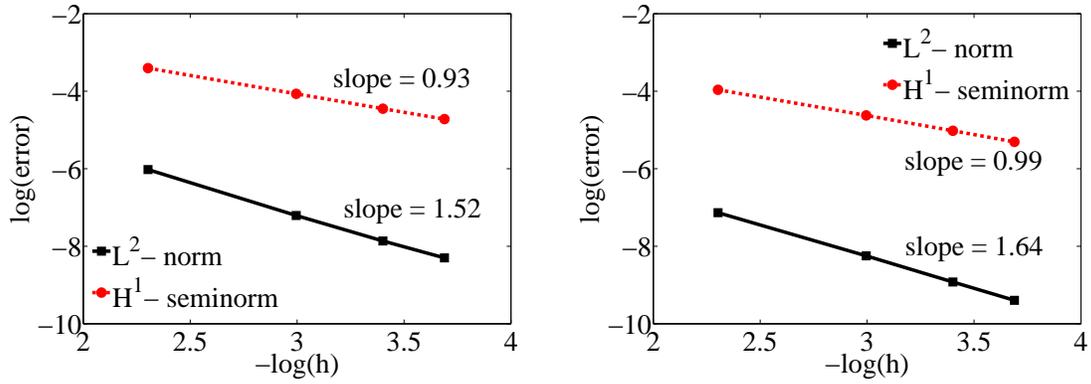

FIGURE 5. Numerical *h*-convergence: Convergence of the proposed computational framework with respect to the *concentration field* is illustrated in this figure. We show the convergence in $L^2$-norm and $H^1$-seminorm for three-node triangular (left) and four-node quadrilateral (right) finite elements.

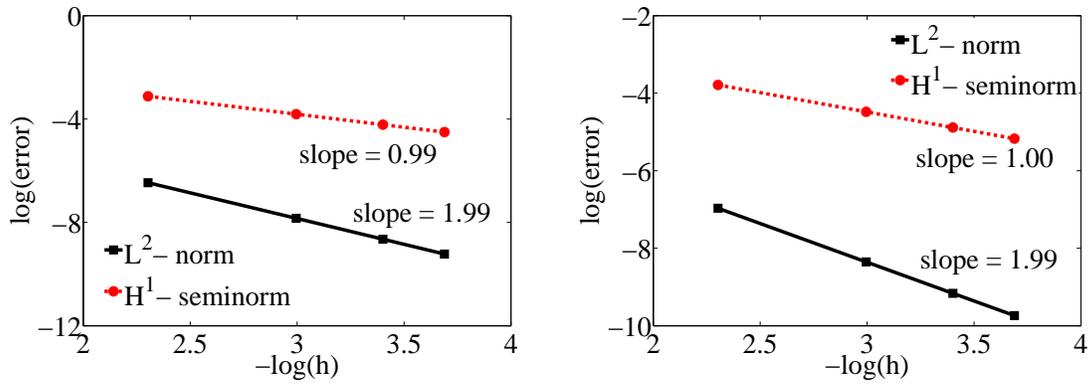

FIGURE 6. Numerical *h*-convergence: Convergence of the proposed computational framework with respect to the *displacement field* is illustrated in this figure. We show the convergence in $L^2$-norm and $H^1$-seminorm for three-node triangular (left) and four-node quadrilateral (right) finite elements.

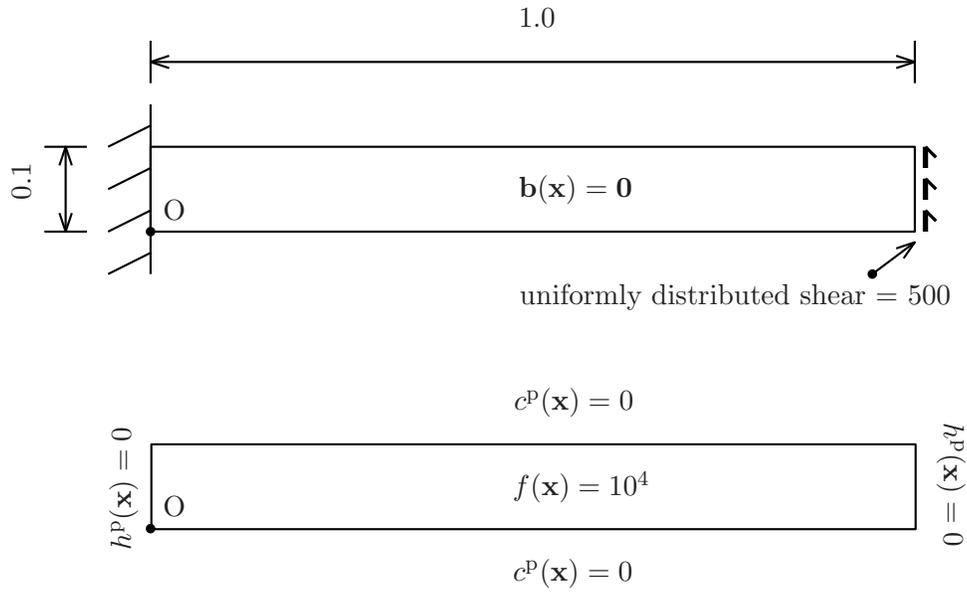

FIGURE 7. Cantilever beam with edge shear: A pictorial description of the boundary value problem with origin at 'O'. The boundary conditions and the specific body force for the deformation problem is shown in the top figure and the boundary conditions and the volumetric source for the diffusion problem is shown in the bottom figure.

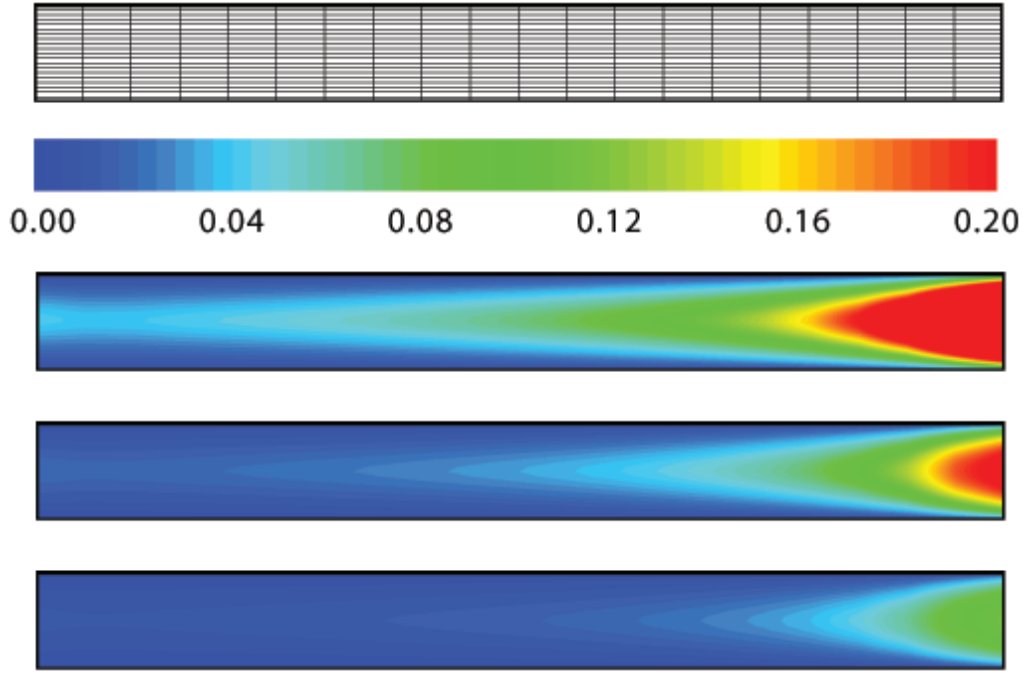

FIGURE 8. Cantilever beam with edge shear: Four-node quadrilateral mesh (top figure) used in the numerical study and comparison of the concentration profiles (in the order of precedence) for three different values of $\Phi_S$ being equal to 5, 10 and 20.

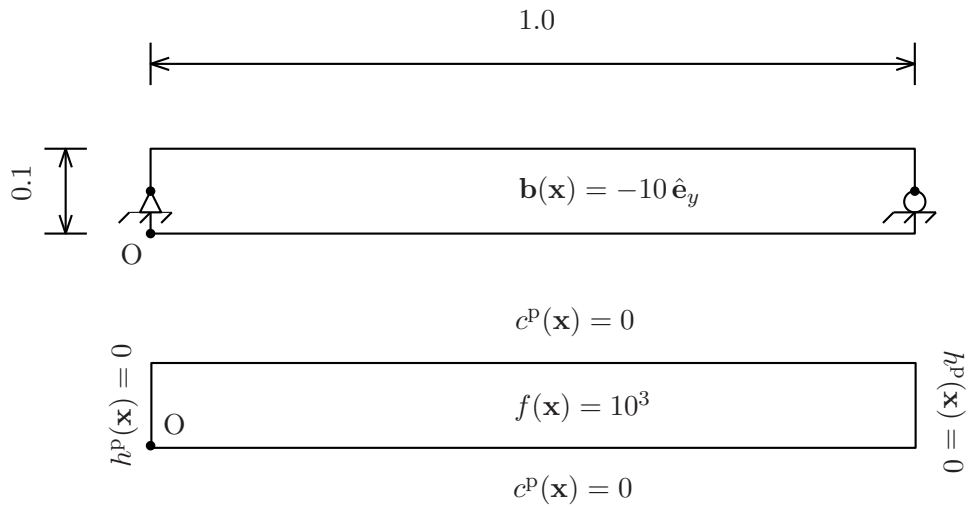

FIGURE 9. Simply supported beam under self-weight: A pictorial description of the boundary value problem with origin at 'O'. The boundary conditions and the specific body force for the deformation problem is shown in the top figure and the boundary conditions and the volumetric source for the diffusion problem is shown in the bottom figure.

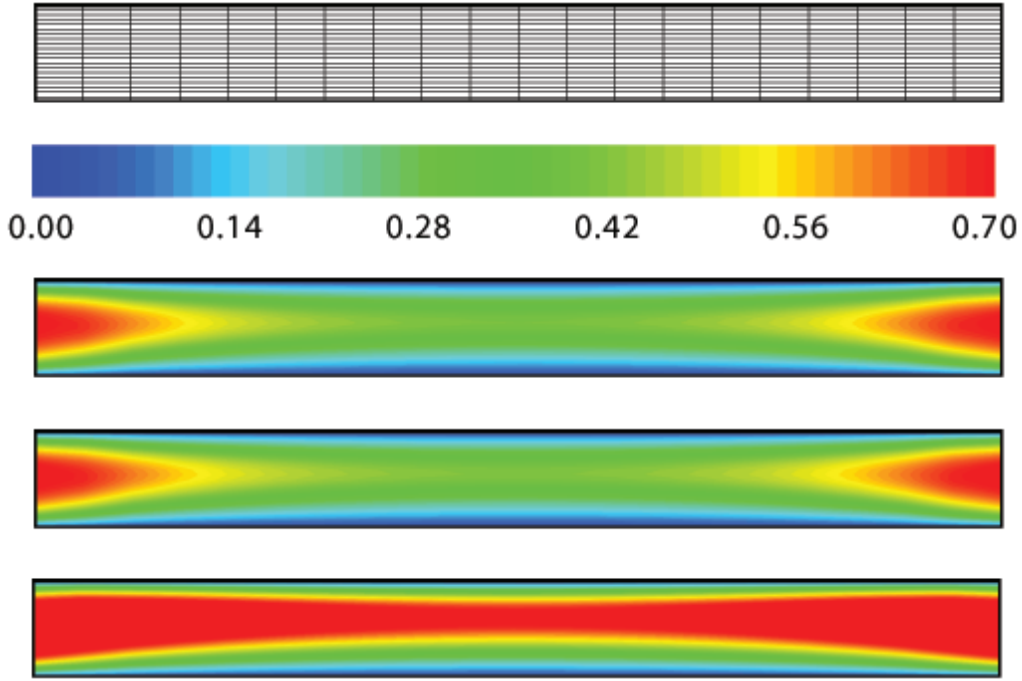

FIGURE 10. Simply supported beam under self-weight: Four-node quadrilateral mesh (top figure) used in the numerical study and comparison of the concentration profiles (in the order of precedence) for three different values of $\eta_S$ being equal to 1, $10^3$ and $2 \times 10^4$.

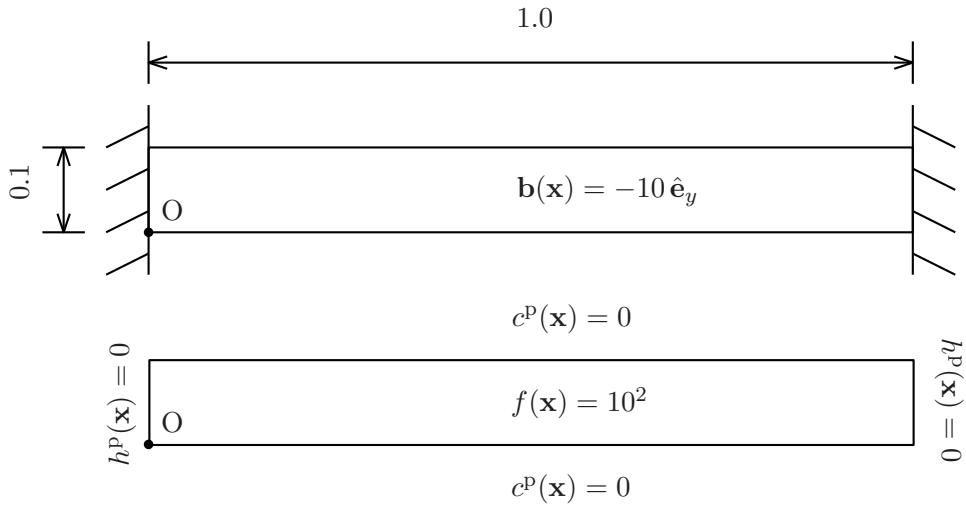

FIGURE 11. Fixed beam under self-weight: A pictorial description of the boundary value problem with origin at 'O'. The boundary conditions and the specific body force for the deformation problem is shown in the top figure and the boundary conditions and the volumetric source for the diffusion problem is shown in the bottom figure.

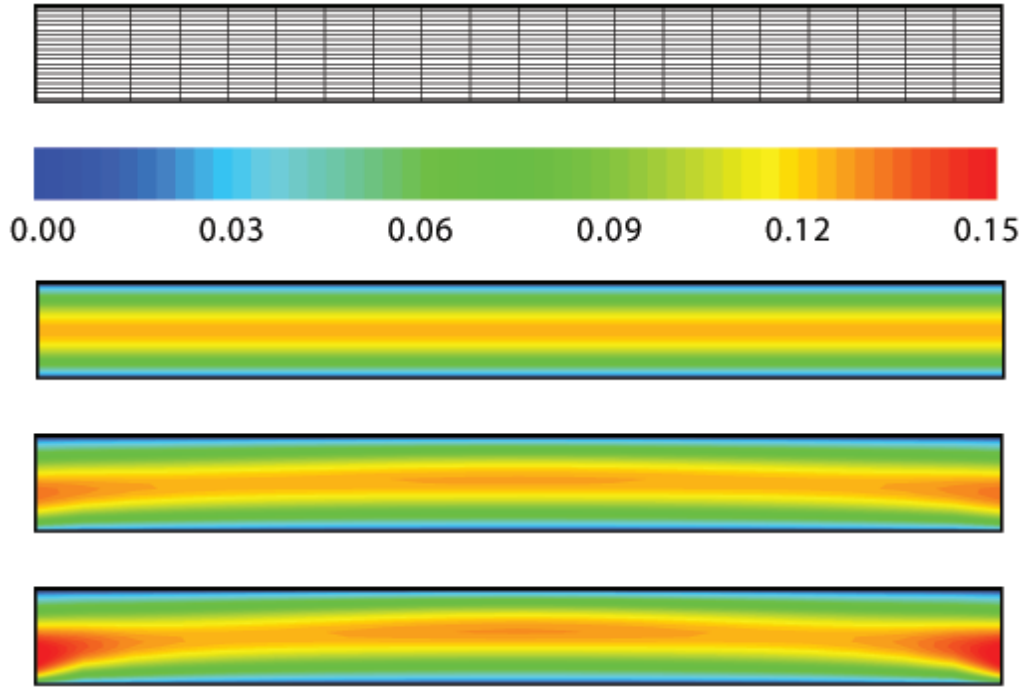

FIGURE 12. Fixed beam under self-weight: Four-node quadrilateral mesh (top figure) used in the numerical study and comparison of the concentration profiles (in the order of precedence) for three different values of $\Phi_T$ being equal to 1, 5 and 7.

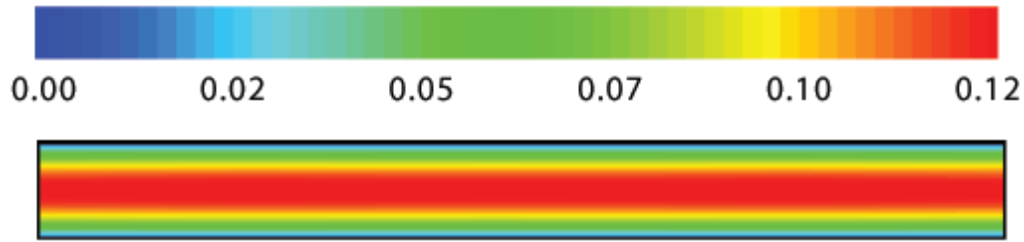

FIGURE 13. Fixed beam under self-weight: The concentration profiles obtained using the diffusivity model (45) as outlined in Reference [25].

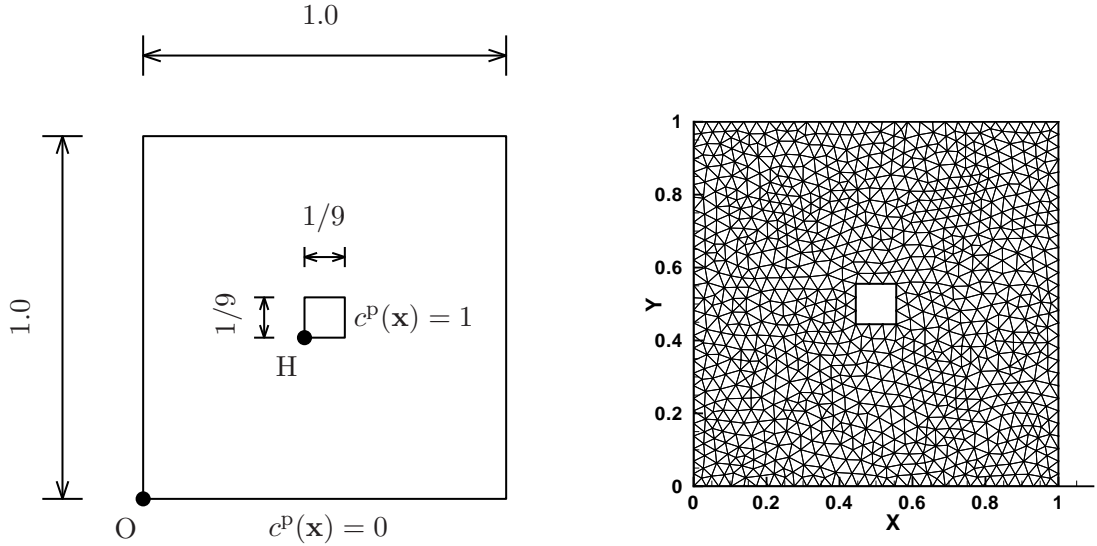

FIGURE 14. Plate with a hole under self-weight: A pictorial description of the dimensions for a plate with a hole in the left figure and three-node triangular unstructured mesh used in our numerical study in right figure. The origin is located at 'O' and the vertex 'H' of the hole at $\left(\frac{4}{9}, \frac{4}{9}\right)$. Displacements at the boundary of the hole and traction at the boundary of the plate are zero. In the region between the plate and the hole, the body force is $\mathbf{b}(\mathbf{x}) = -10\,\hat{\mathbf{e}}_y$ and the volumetric source is $f(\mathbf{x}) = 0$. The concentration at the boundary of the hole is 1 and at the boundary of the plate is 0.

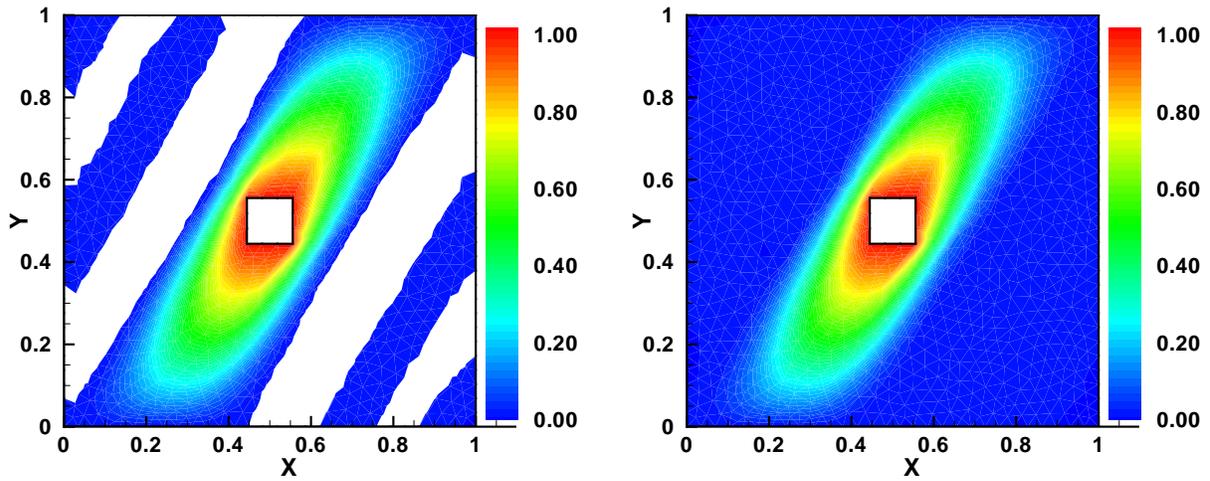

FIGURE 15. Plate with a hole under self-weight: This figure compares concentration profiles from the standard Galerkin formulation (left figure) and the non-negative formulation (right figure). The white area in the left figure indicates the regions in which concentration is negative. The violations also include that are of order machine precision.

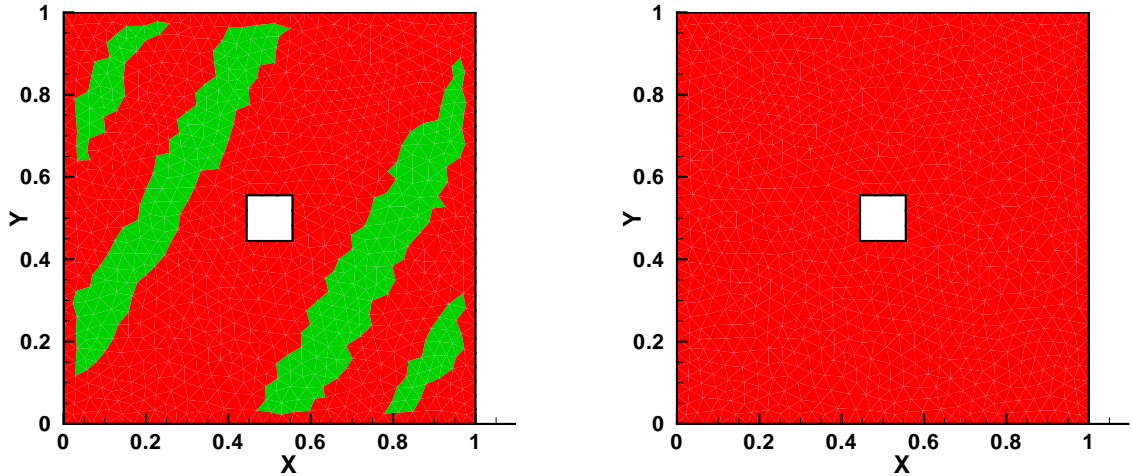

FIGURE 16. Plate with a hole under self-weight: This figure compares degradation profiles from the standard Galerkin formulation (left figure) and the non-negative formulation (right figure). The red area indicates the regions in which the material is degrading and the green area indicates the regions in which the material is healing (which is because unphysical negative value for the concentration). In this figure we do not consider negative values of order machine precision as violations.

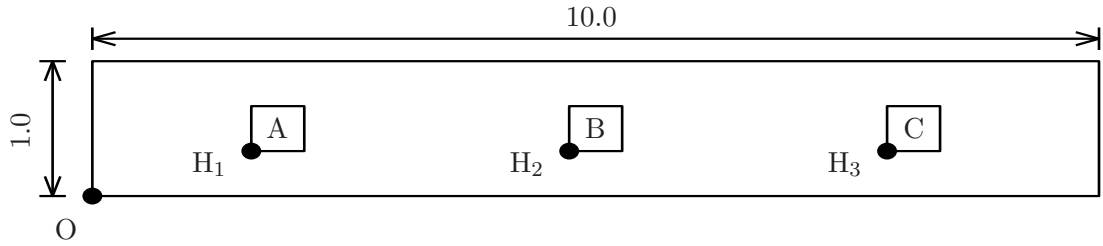

FIGURE 17. Beam with three holes under self-weight: A pictorial description of the dimensions for a beam with three square holes. The dimensions of each hole are $0.4 \times 0.4$. The origin is located at 'O' and the vertices '$H_1$', '$H_2$' and '$H_3$' of the holes 'A, B, C' are at (1.8, 0.3), (4.8, 0.3) and (7.8, 0.3). Displacements at the boundary of the holes and traction at the boundary of the beam are zero. In the region between the beam and the holes, the body force is $\mathbf{b}(\mathbf{x}) = -10\,\hat{\mathbf{e}}_y$ and the volumetric source is $f(\mathbf{x}) = 0$. The concentration at the boundary of the holes is 1 and at the boundary of the beam is 0.

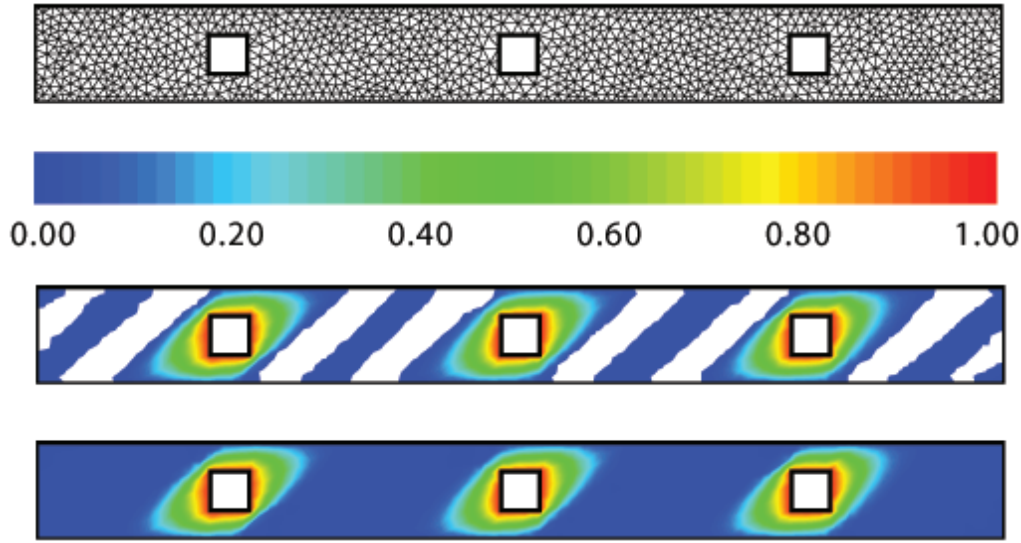

FIGURE 18. Beam with three holes under self-weight: Comparison of concentration profile from standard Galerkin formulation (middle figure) to that of the non-negative formulation (bottom figure). Three-node triangular unstructured mesh used in the numerical study is shown in the top figure. The white area in the middle figure indicates the region in which concentration has a negative value. In this figure, the white area also includes the regions with concentration of order machine precision.

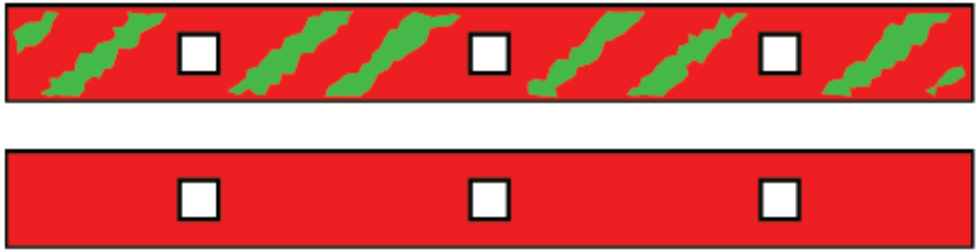

Figure 19. Beam with three holes under self-weight: Comparison of degradation profile from standard Galerkin formulation (top figure) to that of the non-negative formulation (bottom figure). The red area indicates material is degrading and the green area indicates the material is healing (because of unphysical negative value for the concentration). In this figure we do not consider the negative values that are of order machine precision as violations.

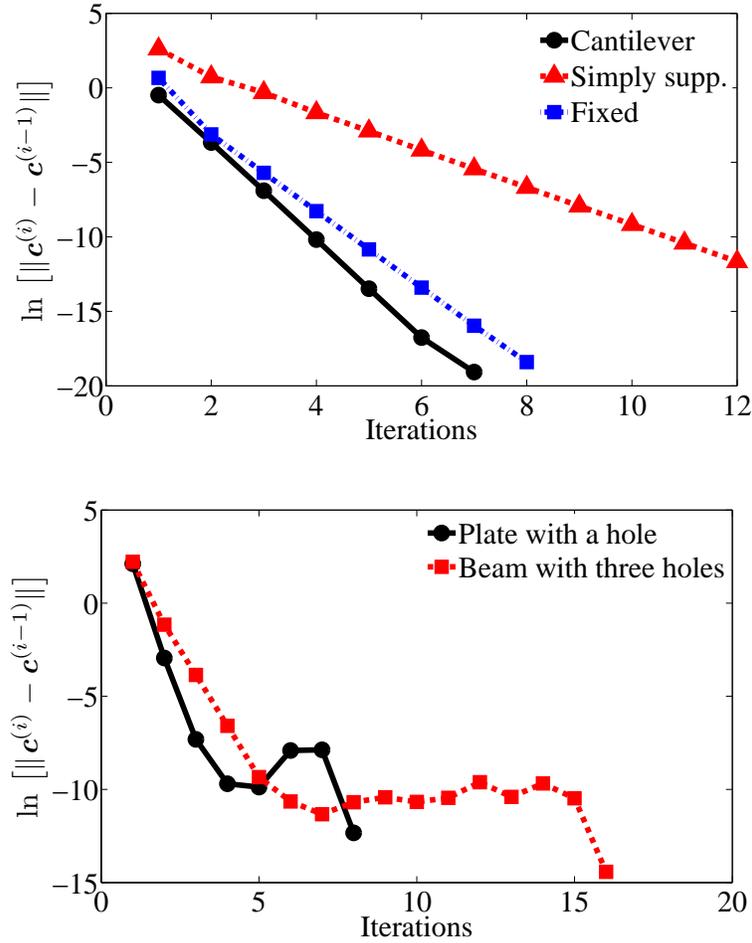

Figure 20. Convergence of staggered coupling algorithm: We have plotted $\ln[\|\boldsymbol{c}^{(i)} - \boldsymbol{c}^{(i-1)}\|]$ with respect to iteration number for various problems. The stopping criterion in the staggered coupling algorithm is $\|\boldsymbol{c}^{(i)} - \boldsymbol{c}^{(i-1)}\| < \epsilon_{\mathrm{TOL}}^{(c)}$. As one can see from the figure, the staggered coupling algorithm converges for all the chosen problems. Also, note that the algorithm need not converge monotonically, which is evident in the case of plate with a hole test problem.

Maruti Kumar Mudunuru, Graduate Student, Department of Mechanical Engineering, Texas A&M University, College Station, Texas 77843.

*E-mail address*: `maruti.iitm@neo.tamu.edu`

Correspondence to: Dr. Kalyana Babu Nakshatrala, Department of Mechanical Engineering, 216 Engineering/Physics Building, Texas A&M University, College Station, Texas 77843. TEL:+1-979-845-1292

*E-mail address*: `knakshatrala@tamu.edu`

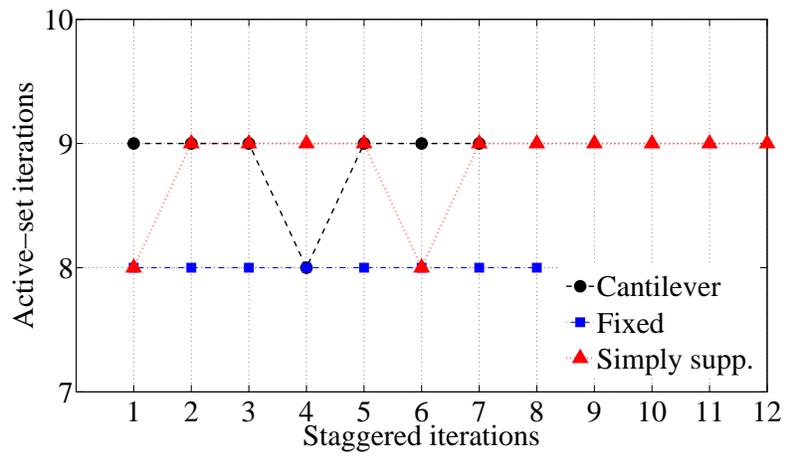

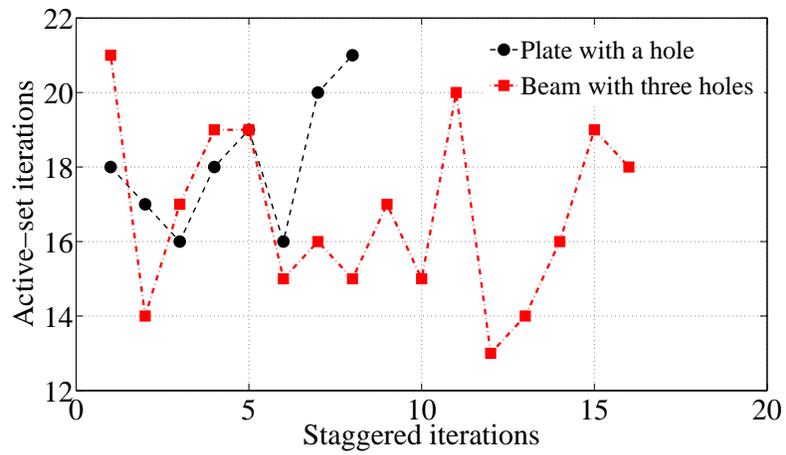

Figure 21. Convergence of the active-set strategy: This figure shows the number of iterations taken by the active-set strategy for each iteration in the staggered coupling algorithm for various test problems.